\documentclass[final]{elsarticle}
\usepackage{lineno,hyperref}
\journal{}

\usepackage{amsthm,amssymb,epsfig,float,color}
\usepackage{siunitx}
\usepackage{upgreek}
\usepackage{geometry}
\newtheorem{lemma}{Lemma}
\newtheorem{theorem}{Theorem}

\usepackage{lipsum}
\usepackage{amsfonts}
\usepackage{booktabs} 
\usepackage{graphicx}
\usepackage{subfigure}
\graphicspath{{Graphics/}}
\usepackage{epstopdf}
\usepackage{algorithmic,multirow,mathtools}

\usepackage{float}
\usepackage{hyperref}
\usepackage{color}
\usepackage{tcolorbox}
\tcbuselibrary{xparse}
\usepackage{relsize}
\usepackage{latexsym}
\usepackage{import}
\usepackage{xifthen}
\usepackage{transparent}

\allowdisplaybreaks

\newcommand{\rmd}{\mathrm{d}}

\newcommand{\pp}[2]{\fracd{\partial {#1}}{\partial {#2}}}
\newcommand{\dd}[2]{\fracd{\rmd {#1}}{\rmd {#2}}}
\newcommand{\fracd}[2]{\displaystyle
{\frac{{\displaystyle{#1}}}{{\displaystyle{#2}}}}}

\newcommand{\calH}{\mathcal{H}}

\newcommand{\Qe}{Q_{\mathrm{e}}}

\newcommand{\Qf}{Q_{\mathrm{f}}}

\newcommand{\Qu}{Q_{\mathrm{u}}}

\newcommand{\Xc}{X_{\mathrm{c}}}
\newcommand{\Xe}{X_{\mathrm{e}}}
\newcommand{\Xf}{X_{\mathrm{f}}}

\newcommand{\bC}{\boldsymbol{C}}

\newcommand{\bPhi}{\boldsymbol{\Phi}}
\newcommand{\bPsi}{\boldsymbol{\Psi}}

\newcommand{\bR}{\boldsymbol{R}}
\newcommand{\bSe}{\boldsymbol{S}_{\mathrm{e}}}
\newcommand{\bSf}{\boldsymbol{S}_{\mathrm{f}}}

\newcommand{\bS}{\boldsymbol{S}}
\newcommand{\bU}{\boldsymbol{\mathcal{U}}}

\newcommand{\bzero}{\boldsymbol{0}}

\newcommand{\qjph}{q_{j+1/2}}

\newcommand{\sme}{\sigma_{\mathrm{e}}}

\newcommand{\vhs}{v_{\mathrm{hs}}}

\newcommand{\br}{\boldsymbol{r}}
\newcommand{\bsigmaC}{\boldsymbol{\sigma}_\mathrm{\!\bC}}
\newcommand{\bsigmaS}{\boldsymbol{\sigma}_\mathrm{\!\bS}}

\definecolor{ros}{RGB}{148,35,9}   

\newcommand{\bCf}{\bC_{\rm f}}
\newcommand{\bCe}{\bC_{\rm e}}

\DeclareMathAlphabet{\mathbcal}{OMS}{cmsy}{b}{n}

\begin{document}

\begin{frontmatter}

\title{A moving-boundary model of reactive settling in wastewater treatment. Part~2: Numerical scheme}

%\author{Raimund B\"urger\fnref{rbfootnote}}
%\address{CI${}^{\mathrm{2}}$MA and Departamento de Ingenier\'{\i}a Matem\'{a}tica, Facultad de Ciencias F\'{i}sicas y Matem\'{a}ticas, Universidad de Concepci\'{o}n, Casilla 160-C, Concepci\'{o}n, Chile}
%\fntext[rbfootnote]{Since 1880.}

\author[RBaddress]{Raimund B\"urger}
\author[JCaddress]{Julio Careaga}
\author[JCaddress]{Stefan Diehl\corref{mycorrespondingauthor}}
\cortext[mycorrespondingauthor]{Corresponding author, \texttt{stefan.diehl@math.lth.se}}
\author[RBaddress]{Romel Pineda}

\address[RBaddress]{CI${}^{\mathrm{2}}$MA and Departamento de Ingenier\'{\i}a Matem\'{a}tica, Facultad de Ciencias F\'{i}sicas y Matem\'{a}ticas, Universidad~de~Concepci\'{o}n, Casilla 160-C, Concepci\'{o}n, Chile}
\address[JCaddress]{Centre for Mathematical Sciences, Lund University, P.O.\ Box 118, S-221 00 Lund, Sweden}

\begin{abstract}
A numerical scheme is proposed for the simulation  of reactive settling  in sequencing batch reactors (SBRs) in wastewater treatment plants.
Reactive settling is the process of sedimentation of flocculated particles (biomass; activated sludge) consisting of several material components that react with substrates dissolved  in the fluid.
An SBR is operated in cycles of consecutive fill, react, settle, draw and idle stages, which means that the volume in the tank varies and the surface moves with time.
The process is modelled by a system of spatially one-dimensional, nonlinear, strongly degenerate parabolic convection-diffusion-reaction equations.
This system is coupled via conditions of mass conservation to transport equations on a half line whose origin is located at a moving boundary and that models the effluent pipe.
A finite-difference scheme  is proved to satisfy an invariant-region property (in particular, it is positivity preserving) if executed in a simple splitting way. 
Simulations are presented with a modified variant of the established activated sludge model no.~1 (ASM1).
\end{abstract}


\end{frontmatter}

%\linenumbers

\section{Introduction}

\begin{figure}[t]
\centering 
\includegraphics[scale=1.1]{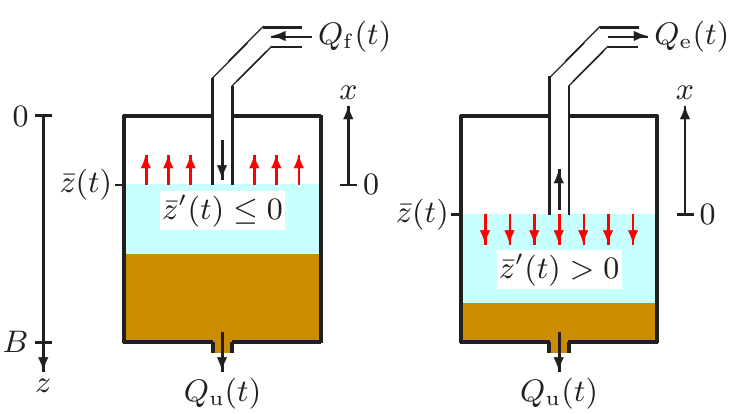}
\caption{Fill at the volume rate $\Qf(t)>0$ greater than the underflow rate $\Qu(t)\geq 0$ resulting in a rise of the mixture surface location~$z=\bar{z}(t)$. Right: Draw (extraction) of mixture from the surface at the rate $\Qe(t)>0$ implies a descending surface.}\label{fig:FillDraw}
\end{figure}

\begin{figure}[t]
\centering
 \includegraphics[width=0.99\textwidth]{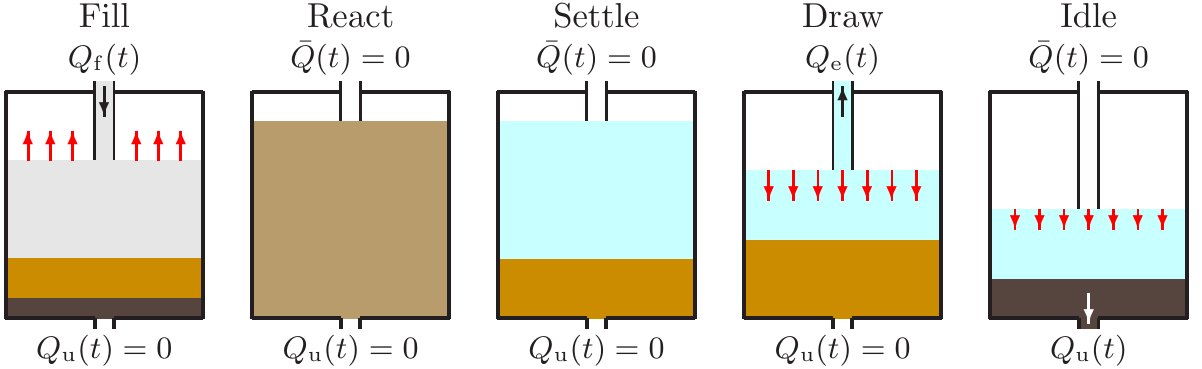}
\caption{The five stages of a cycle of an SBR.
The tank is first filled with wastewater at the volumetric flow $\Qf(t)>0$ and concentrations~$\bCf(t)$ and $\bSf(t)$.
During the react stage, biological reactions take place under complete mixing by an impeller or by aeration.
Then batch sedimentation with reactions occurs and liquid is extracted during the draw stage.
During the idle stage, some of the bottom sludge can be withdrawn and then the fill stage starts again.}
\label{fig:SBRcycle}
\end{figure}

A sequencing batch reactor (SBR) is a tank (with possibly varying cross-sectional area) used for the purification of wastewater (e.g., \cite{chen2020book,droste,metcalf}).
It has a controlled outlet at the bottom and at the surface of the mixture, a floating device allows for controlled fill or extraction of mixture; see Figure~\ref{fig:FillDraw}.
It is operated batch-wise in a sequence of cycles of fill, react, settle, draw and idle stages; see see Figure~\ref{fig:SBRcycle}.
The mixture consists of flocculated particles of biomass (activated sludge) consisting of several components that react with dissolved substrates (nutrients) in the liquid. (The removal of these substrates is the purpose of an SBR.) 
The full model of the process consists of a system of nonlinear partial differential equations (PDEs); see~\cite{bcdp_part1}, to which we also refer for references to previous related works.
To the authors' knowledge, models for SBRs with reactions during all stages (including the settle) 
 have rarely been proposed in the literature.  
Models of reactive settling in continuously operated secondary settling tanks (SSTs) based on PDEs and numerical schemes are presented in~\cite{bcdimamat21,SDcace_reactive,SDm2an_reactive}; see also references therein.

In addition to the mathematical and numerical difficulties of an SST model, the moving boundary in an SBR means a special challenge.
Although in normal operation only liquid is extracted through the pipe in the draw phase, our model and numerical scheme are capable to handle any concentrations at the mixture surface.

The moving-boundary model has a connected half-axis with transport equations and nonlinear mass-preserving coupling conditions that    do not define the coupling concentrations uniquely.
Such a problem of nonuniqueness arises already for a scalar conservation law with discontinuous flux, which has been investigated widely~\cite{Adimurthi&J&G2004,Andreianov2011,SDsiam1,Gimse1992,Karlsen&R&T2002,Karlsen&R&T2003}; in particular, in the context of  continuous sedimentation~\cite{Burger&K&T2005a,SDcec_varyingA,Chancelier1994,SDsiam2,SDsiam3} where a monotone numerical scheme approximates the correct solutions~\cite{Burger&K&T2005a}.

It is the purpose of this contribution to present a numerical scheme that can handle  the strong type degeneracy of the PDEs as well as  the moving boundary where both a source is located and a half-line model attached.
The scheme is monotone when the reaction terms are zero.
We prove that if the scheme is used in the simple Lie-Trotter-Kato splitting way (e.g.\ \cite{Holden_splitting2010}), namely, one explicit time step is taken without reactions and a another step with only reactions, then the numerical solutions have an invariant-region property under a convenient Courant-Friedrichs-Lewy (CFL) condition.
In particular, concentrations are nonnegative.

The PDE model is presented in Section~\ref{sec:model}.
Section~\ref{sec:numscheme} contains the numerical scheme based on the one proposed  in~\cite{bcdimamat21} but modified to handle the moving boundary. 
The fully discrete and explicit scheme is presented in Section~\ref{subsec:scheme} and its splitting version with the CFL condition and invariant-region property is presented in Section~\ref{subsec:numcfl}. 
The numerical scheme during the full mixing stage can be found in Section~\ref{sec:numfullmix}.
In Section~\ref{sec:numex}, we show a numerical examples of SBR operation with a constant cross-sectional area (cylindrical vessel)
and a commonly used activated sludge model for the biokinetic reactions; however, slightly adjusted to deliver non-negative concentrations only.
Some conclusions can be found in Section~\ref{sec:concl}.

\section{The model}\label{sec:model}

The vertical  $z$-axis of the governing model and the moving coordinate system with the $x$-half-axis are  shown in Figure~\ref{fig:FillDraw}.
The characteristic function~$\gamma$ equals one inside the mixture and zero otherwise, i.e., $\gamma(z,t)=\chi_{\{\bar{z}(t)<z<B\}}$, where $\chi_I$ is the indicator function which equals one if and only if $I$ is true, and $z=\bar{z}(t)$ is the surface location.
We let $A=A(z)$ denote the cross-sectional area of the tank that may depend on depth~$z$.

The solid phase consists of flocculated particles of $k_{\boldsymbol{C}}$ types   with concentrations~$\smash{C^{(1)} ,\dots, C^{(k_{\boldsymbol{C}})}}$.
 The components of the liquid phase are water of concentration~$W$ and $k_{\bS}$ dissolved substrates of concentrations $\smash{S^{(1)} ,\dots, S^{(k_{\bS})}}$. 
The total concentrations of solids~$X$ and liquid~$L$ are
\begin{align} \label{eq:XLdef} 
 X :=  C^{(1)} + \dots  + C^{(k_{\boldsymbol{C}})} ,\quad L := 
  W + S^{(1)} + \dots + S^{(k_{\boldsymbol{S}})}.
\end{align} 
All these concentrations depend on~$z$ and~$t$.

At the surface of the mixture, $z=\bar{z}(t)$, we model a floating device connected to a pipe through which one can feed the tank with a given volume rate~$\Qf(t)$ and given feed concentrations $\bCf(t)$ and $\bSf(t)$; see Figure~\ref{fig:FillDraw}.
Alternatively, this floating device allows to extract mixture at a given volume rate $\Qe(t)>0$ through the same pipe; hence, one cannot fill and extract simultaneously.
If~$[0,T]$ denotes the total time interval of modelling (and simulation in Section~\ref{sec:numscheme}), we assume that $T:={T}_\mathrm{e}\cup{T}_\mathrm{f}$, where
\begin{align*}
{T}_\mathrm{e}&:= \bigl\{t\in\mathbb{R}_+:\Qe(t)>0,\Qf(t)=0 \bigr\},
\qquad
{T}_\mathrm{f}:= \bigl\{t\in\mathbb{R}_+:\Qe(t)=0,\Qf(t)\geq0 \bigr\}.
\end{align*}
When~$t\in{T}_\mathrm{e}$, we model in~\cite{bcdp_part1} the extraction flow in the effluent pipe by a moving coordinate system, namely  a half line~$x\geq 0$, where $x=0$ is attached to $z=\bar{z}(t)$; see Figure~\ref{fig:FillDraw}.
Since we assume that there are no reactions in the pipe and all components have the same velocity, the conservation law for the pipe is the linear advection equation $A_\mathrm{e}\partial_t\boldsymbol{\tilde{C}}+\Qe(t)\partial_x\boldsymbol{\tilde{C}}=0$ where $A_\mathrm{e}$ is the  cross-sectional area of the pipe. However, we  
 are only interested in the effluent concentrations in the pipe at $x=0^+$, which we denote by~$\bC_\mathrm{e}(t)$ and $\bS_\mathrm{e}(t)$.
At the bottom, $z=B$, one can withdraw mixture at a given volume rate $\Qu(t)\geq 0$.
The underflow region $z>B$ is for simplicity modelled by setting $A(z):= A(B)$, since we are only interested in the underflow concentration $\bC_\mathrm{u}(t)$, which is an outcome of the model (analogously for~$\bS_\mathrm{u}(t)$).

Below the surface, the volume fractions of solids and liquid add to one, i.e.,
\begin{equation*}
\frac{X}{\rho_X} + \frac{L}{\rho_L} = 1
\quad\Leftrightarrow\quad
L=\rho_L\left(1-\frac{X}{\rho_X}\right).
\end{equation*}
The same holds for the feed concentrations.
For known   $\bC$ and $\bS$, \eqref{eq:XLdef} implies   the water concentration  
\begin{equation}  \label{eq:W}
W = \rho_L(1-X/\rho_X)-\big(S^{(1)}+\cdots+S^{(k_{\bS})}\big).
\end{equation} 
This concentration is not part of any reaction and can be computed afterwards.

For computational purposes, we define a maximal concentration~$\hat{X}$ of solids and  assume that the density of all solids is the same, namely~$\rho_X>\hat{X}$.
Similarly, we assume that the liquid phase has the density~$\rho_L<\rho_X$, typically the density of water.
The velocity of a solid particle $v_X = q+v$ is the sum of the average volumetric velocity, or bulk velocity, of the mixture
\begin{equation} \label{eq:qdef2}
q(z,t) :=\frac{\Qu(t)}{A(z)}\chi_{\{z>\bar{z}(t)\}},
\end{equation}
and the excess velocity~$v$, which is given by the following commonly used expression~\cite{Burger&K&T2005a,SDwatres3}:
\begin{equation}\label{eq:v}
v:=v(X,\partial_z X,z,t) 
:= \gamma(z,t)\vhs(X)\left(
1-\dfrac{\rho_X\sme'(X)}{Xg\Delta\rho}\partial_z{X}
\right)
= \gamma(z,t)\big(\vhs(X) - \partial_z{D(X)}\big),
\end{equation}
where
\begin{equation*}
D(X) := \int_{X_c}^{X}d(s)\,{\rm d}s,
\qquad
d(X):=\vhs(X)\dfrac{\rho_X\sme'(X)}{gX\Delta\rho}.
\end{equation*}
Here, $\Delta\rho:=\rho_X-\rho_L$, $g$ is the acceleration of gravity, $\vhs=\vhs(X)$ is the hindered-settling velocity, which is assumed to be decreasing and satisfy $\vhs(\hat{X})=0$,
$\sme=\sme(X)$ the effective solids stress, which satisfies $\sme'(X) = 0$ for $X \leq \Xc$ and $\sme'(X) > 0$ for $X>\Xc$,  
where $\Xc$ is a critical concentration above which the particles touch each other and form a network that can bear a certain stress.

The reaction terms for all components are collected in the vectors 
\begin{equation}\label{eq:reactionrates}
\bR_{\bC}(\bC,\bS)=\bsigmaC\br(\bC,\bS),\qquad
\bR_{\bS}(\bC,\bS)=\bsigmaS\br(\bC,\bS),
\end{equation}
which model the increase of solid and soluble components, respectively, where $\bsigmaC$ and $\bsigmaS$ are constant stoichiometric matrices and $\br(\bC,\bS)\geq\bzero$ is a vector of non-negative reaction rates, which are assumed to be bounded and Lipschitz continuous functions.
We set 
\begin{align*} 
\tilde{R}_{\bC}(\bC,\bS) :=  R_{\bC}^{(1)}(\bC,\bS)
+\dots + R_{\bC}^{(k_{\bC})}(\bC,\bS)
\end{align*} 
 (analogously for~$\smash{\tilde{R}_{\bS}(\bC,\bS)}$).
 
 In the present work the matrices~$\bsigmaC$ and~$\bsigmaS$ and the vector~$\boldsymbol{r}$ are given by expressions that represent a  slightly modified version of the activated sludge model no.~1, usually abbreviated `ASM1' in wastewater engineering. The model was the first comprehensive activated sludge model developed by a task group of the International Water Association (IWA) (see \cite{Henze1987WR}). The model describes the principal biochemical reactions that occur within the activated sludge process, 
  which is the most widespread technology for the secondary treatment of municipal wastewater and as commented in~\cite{makinia20}, 
   constitutes `the heart' of many wastewater treatment plants. Later versions, known as ASM2, ASM2d, and ASM3, account for 
    additional reactions such as fermentation and chemical or biological phosphorus removal (not considered in the present work). 
     Commercial software packages that include these models are used commonly by wastewater process design engineers 
      for the process design of various activated sludge system configurations \cite{metcalf}.
We refer to Appendix~A for the description of the modified ASM1 model used herein and to \cite{chen2020book,metcalf,makinia20} for further information on activated sludge models.

In order to establish an invariant-region property for the numerical solution, we make some technical assumptions.
To ensure that  the numerical solution for the solids does not exceed the maximal concentration~$\hat{X}$, we assume the following:
\begin{align} 
&\text{there exists an $\varepsilon>0$ such that $\bR_{\bC}(\bC,\boldsymbol{S})=\bzero$ for all $X\geq\hat{X}-\varepsilon$},\label{eq:techRC}\\
&v(\hat{X},\partial_z X,z,t)=0.\label{eq:techrel}
\end{align} 
These conditions mean that when the concentration is (near) the maximal one, biomass cannot grow any more and its relative velocity to the liquid phase is zero.
To obtain positivity of component~$k$ of the concentration vector~$\bC$, we let
\begin{equation*}
I_{\bC,k}^-:=\big\{l\in\mathbb{N}:\sigma_{\bC}^{(k,l)}<0\big\},
\qquad
I_{\bC,k}^+:=\big\{l\in\mathbb{N}:\sigma_{\bC}^{(k,l)}>0\big\},
\end{equation*}
denote the sets of indices~$l$ that have negative and positive stoichiometric coefficients, respectively, and assume the following (analogously for~$\bS$):
\begin{equation}\label{eq:positivity}
\text{if $l\in I_{\bC,k}^-$, then $r^{(l)}(\bC,\bS)=\bar{r}^{(l)}(\bC,\bS)C^{(k)}$ with $\bar{r}^{(l)}$ bounded.}
\end{equation}
Assumption~\eqref{eq:positivity} implies that 
\begin{equation*} 
R^{(k)}_{\bC}(\bC,\bS)\bigr|_{C^{(k)}=0} \geq 0 \quad \text{for $k=1, \dots, k_{\boldsymbol{C}}$}\label{eq:assumptionRC}
\end{equation*}
(analogously for~$\bS$), which means that the system of ODEs 
\begin{equation}\label{eq:systODE}
\frac{\rmd}{\rmd t} \begin{pmatrix}
\bC\\ \bS
\end{pmatrix} = \begin{pmatrix}
\bR_{\bC}(\bC,\bS)\\ \bR_{\bS}(\bC,\bS)
\end{pmatrix}
\end{equation}
has a non-negative solution if the initial data are non-negative~\cite{Formaggia2011}.
This positivity property is carried over to the numerical splitting scheme suggested here.

The volume of the mixture is defined by 
\begin{equation*}
\bar{V}(t):= V \bigl(\bar{z}(t)\bigr), \quad \text{where} \quad 
V(z):= \int_{z}^{B}A(\xi)\,\mathrm{d}\xi \quad \text{for $ 0\leq z\leq B$.} 
\end{equation*} 
The function~$V$  is invertible since $V'(z)=-A(z)<0$; in particular,
\begin{equation} \label{eq:Vprime}
\bar{V}'(t)=V' \bigl(\bar{z}(t) \bigr)\bar{z}'(t)=-A \bigl(\bar{z}(t) \bigr)\bar{z}'(t).
\end{equation} 
It turns out that the surface location is given by the following explicit expression of the given volumetric flows:
\begin{equation*}
\bar{z}(t)=V^{-1}\left(
\bar{V}(0)+ \int_{0}^{t}\big(\bar{Q}(s) -\Qu(s)\big)\,\mathrm{d}s\right),
\qquad
\bar{Q}(t):=\begin{cases}
-\Qe(t)<0 & \text{if $t\in {T}_\mathrm{e}$,} \\
\Qf(t)\geq 0 & \text{if $t\in {T}_\mathrm{f}$.}
\end{cases}
\end{equation*} 
Alternatively, $\bar{z}(t)$ can be obtained from
\begin{equation} \label{eq:zbarprime}
\bar{z}'(t)=\frac{\Qu(t)-\bar{Q}(t)}{A(\bar{z}(t))}.
\end{equation}

To state the governing model, we define the velocities
\begin{equation} \label{eq:modelterms} 
\begin{split} 
 \mathcal{F}_{\boldsymbol{C}} (X,z,t)  & :=  q(z,t) + \gamma(z,t) \vhs(X),  \\
 \mathcal{F}_{\boldsymbol{S}} (X,z,t) & := \frac{ \rho_X q(z,t) - ( q(z,t) + \gamma(z,t) \vhs(X) )  X }{\rho_X - X},
\end{split}
\end{equation}
and then introduce the total mass fluxes as
\begin{align}
\bPhi_{\bC}  := \bPhi_{\bC}(\bC,X,\partial_z X,z,t) &:= A(z)\big(\mathcal{F}_{\boldsymbol{C}}(X,z,t) - \gamma(z,t)\partial_z{D(X)}\big), \\
\bPhi_{\bS}  := \bPhi_{\bS}(\bS,X,\partial_z X,z,t) &:= A(z)\mathcal{F}_{\boldsymbol{S}}(X,z,t), \\
\bPhi_{\bC,\mathrm{e}}(z,t)
&:= \left.\Big(
A(z)\big(
\vhs(X) - \partial_z{D(X)}\big)
- \Qe
\Big)\bC
\right|_{z=\bar{z}(t)^+},\label{eq:PhiC_e}\\
\bPhi_{\bS,\mathrm{e}}(z,t)
&:= \left.
- \left(A(z)
 \frac{X(\vhs(X) - \partial_z{D(X)})}{\rho_X-X}
+ \Qe\right)\bS
\right|_{z=\bar{z}(t)^+}.\label{eq:PhiS_e}
\end{align}
The complete  model is the following ($\delta$ is the delta function):
\begin{subequations}\label{finalmod}
\begin{alignat}2
A(z)\partial_t{\bC}+\partial_z{\bPhi_{\bC}} & = \delta \bigl(z-\bar{z}(t)\bigr)\Qf\bCf + \gamma(z,t)A(z)\bR_{\bC}, &\quad& z\in\mathbb{R}, \label{finalmod_a} \\
A(z)\partial_t{\bS} +\partial_z{\bPhi_{\bS}} & = \delta \bigl(z-\bar{z}(t)\bigr)\Qf\bSf
+ \gamma(z,t)A(z)\bR_{\bS},&\quad& z\in\mathbb{R}, \label{finalmod_b}\\
\bC_\mathrm{e}(t) &= -{\bPhi_{\bC,\mathrm{e}}(z,t)}/{\Qe(t)},&\quad& t\in{T}_\mathrm{e},\label{eq:connC}\\
\bS_\mathrm{e}(t) &= -{\bPhi_{\bS,\mathrm{e}}(z,t)}/{\Qe(t)},&\quad& t\in{T}_\mathrm{e},\label{eq:connS}\\
\bC_{\rm u}(t)&:=\bC(B^+,t),&\qquad& t>0,\label{eq:Cu_conc_def}\\
\bS_{\rm u}(t)&:=\bS(B^+,t),&\qquad& t>0. \label{eq:Su_conc_def}
\end{alignat} 
\end{subequations}
The water concentration $W$ can always be calculated from~\eqref{eq:W}.
No initial data are needed for the outlet concentrations, but for~$\boldsymbol{C}$ and~$\boldsymbol{S}$, namely 
\begin{equation*}
\bC^0= \big(C^{(1),0},C^{(2),0},\dots, C^{(k_{\bC}),0}\big)^{\mathrm{T}},\qquad
\bS^0= \big(S^{(1),0},S^{(2),0},\dots, S^{(k_{\bS}),0}\big)^{\mathrm{T}}.
\end{equation*}

During the react stage of an SBR (see Figure~\ref{fig:SBRcycle}), full mixing occurs and the system of PDEs~\eqref{finalmod_a} and \eqref{finalmod_b} reduces to the following system of ordinary differential equations (ODEs) 
 for the homogeneous concentrations in $\bar{z}(t)<z<B$:
\begin{subequations}\label{eq:ODEmodel}
\begin{align}
\bar{V}(t)\dd{\bC}{t} & = \big(\Qu(t)-\bar{Q}(t)\big)\bC + 
\Qf(t)\bCf(t) + \bar{V}(t)\bR_{\bC},  \label{eq:ODEmodelC} \\
\bar{V}(t)\dd{\bS}{t} & = \big(\Qu(t)-\bar{Q}(t)\big)\bS + 
\Qf(t)\bSf(t) + \bar{V}(t)\bR_{\bS},  \label{eq:ODEmodelS}
\end{align}
\end{subequations}
where all concentrations depend  on time only since they are averages (below the surface).
As before, $W$ can be obtained afterwards from~\eqref{eq:W}.
In the region $0<z<\bar{z}(t)$ all concentrations are zero.
Because of~\eqref{eq:PhiC_e} and \eqref{eq:connC} we have  $\bC_{\rm u}(t)=\bC(t)$ and $\bC_{\rm e}(t)=\bC(t)\chi_{\{t\in{T}_\mathrm{e}\}}$ (analogously for~$\bS$).

\section{Numerical scheme} \label{sec:numscheme}

\subsection{Spatial discretization and numerical fluxes} \label{subsec:spatdisc}

We divide the tank into $N$ computational cells each having  depth $h = B/N$. 
 Assume that the midpoint   of cell~$j$   has the coordinate~$z_j$, hence, the cell is the interval
$[z_{j-1/2},z_{j+1/2}]$. The top cell~1 is thus  $[z_{1/2},z_{3/2}]=[0,h]$,
and the bottom location is $z=z_{N+1/2}=B$.
To obtain the underflow concentrations, we add one cell below~$z=B$.
To obtain the extraction concentrations, we add one cell $[0,\Delta x]$ of the $x$-coordinate system.
To approximate the cell volumes, we define the average cross-sectional areas
\begin{equation*}
 A_{j-1/2} := \dfrac{1}{h}\int_{z_{j-1}}^{z_j} A(\xi)\,{\rm d} \xi \quad \mbox{and}\quad  A_j := \dfrac{1}{h}\int_{z_{j-1/2}}^{z_{j+1/2}} A(\xi)\,{\rm d }\xi.
\end{equation*}

The unknowns are approximated by functions that are piecewise constant in each cell~$j$, i.e.\ $\smash{C^{(k)}(z,t)\approx C_j^{(k)}(t)}$,
$z\in[z_{j-1/2},z_{j+1/2}]$, which are collected in the vector~$\bC_{j}(t)$.
We define $\bar{j}(t): = \lceil \bar{z}(t)/h \rceil$, which is the smallest integer larger than or equal to $\bar{z}(t)/h$. 
Then the surface $z=\bar{z}(t)$ is located in the surface cell~$\bar{j}(t)$.

We let $\gamma_{j+1/2}(t):=\gamma(z_{j+1/2},t)$ and define the approximate volume-average velocity $\qjph(t):= q(z_{j+1/2},t)$ in accordance with~\eqref{eq:qdef2} via
\begin{align*}
A_{j+1/2}\qjph(t) := \Qu(t)\chi_{\{j+1/2>\bar{j}(t)\}}.
\end{align*}

Using the notation $a^-:=\min\{a,0\}$ and $a^+:=\max\{a,0\}$, we define
\begin{align} \label{eq:numfluxes}
\begin{aligned}
  J_{j+1/2}^{\bC} &= J_{j+1/2}^{\bC}(X_j,X_{j+1}) :=  \big(D(X_{j+1})-D(X_j)\big)/h,\\
v^{X}_{j+1/2}&=v^{X}_{j+1/2}(X_{j},X_{j+1},t)
:=\qjph + \gamma_{j+1/2}\big(\vhs(X_{j+1}) - J_{j+1/2}^{\bC}\big),\\
F^X_{j+1/2}&= F^X_{j+1/2}(X_{j},X_{j+1},t) :=  (v_X X)_{j+1/2} := v^{X,-}_{j+1/2} X_{j+1} + v^{X,+}_{j+1/2} X_{j},\\
\bPhi^{\bC}_{j+1/2}&
:= A_{j+1/2}\bigl(v^{X,-}_{j+1/2}\bC_{j+1} + v^{X,+}_{j+1/2}\bC_{j}\bigr),\\
\bPhi^{\bS}_{j+1/2}&
:= A_{j+1/2}\biggl(\dfrac{(\rho_Xq_{j+1/2}-F^X_{j+1/2})^-}{\rho_X-X_{j+1}}\bS_{j+1} + \dfrac{(\rho_Xq_{j+1/2}-F^X_{j+1/2})^+}{\rho_X-X_{j}}\bS_{j}\biggr).
\end{aligned}
\end{align} 
In particular, we have $\bPhi^{\bC}_{j+1/2}=\bPhi^{\bS}_{j+1/2}=\bzero$ for $j<\bar{j}(t)$.
We  denote by $[\Delta \bPhi]_j:=\bPhi_{j+1/2}-\bPhi_{j-1/2}$  the flux difference associated with cell~$j$.
For the single cell on the $x$-axis, the advective fluxes at $x=\Delta x$ are $\Qe(t)\bCe(t)$ and $\Qe(t)\bSe(t)$.

\subsection{Time discretization and surface fluxes} \label{subsec:timedisc}

\begin{figure}[t]
\centering
 \includegraphics[width=\textwidth]{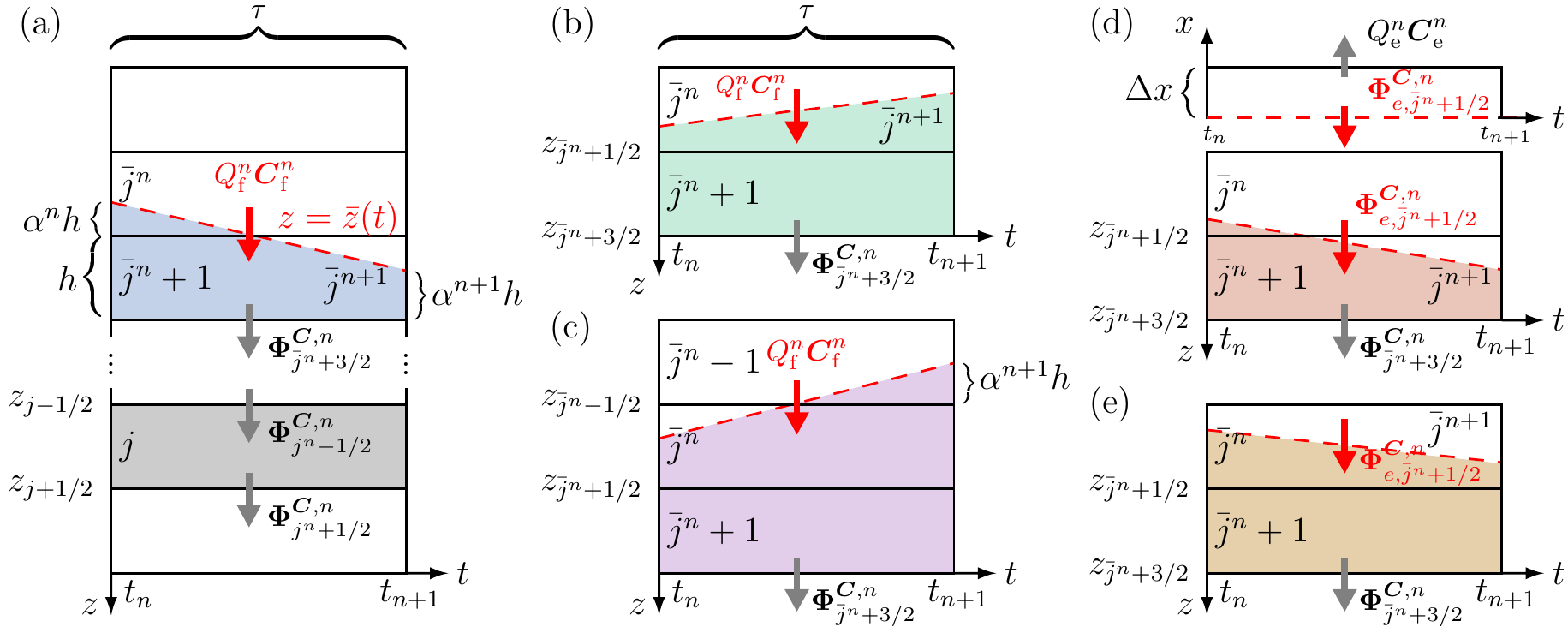}
 \caption{Fluxes over cell boundaries shown by grey arrows and the flux at the surface with red arrows. The surface level~$z=\bar{z}(t)$ is drawn with a red dashed line.
Plot~(d) shows the extraction pipe cell where the origin of the $x$-axis located on the red dashed surface~$z=\bar{z}(t)$ in plots~(d) and~(e).
\label{fig:fig_cells}}
\end{figure}

We let $T$ denote the total simulation time, $t_n$, $n=0,1,\ldots,N_T$, the discrete time points and $ \tau := T/N_T$
the time step that should satisfy a certain CFL condition; see below.
The value of a variable at time $t_n$ is denoted by an upper index, e.g., $\bC_j^n$ and it is thus assumed to be constant in time during $t_n\leq t< t_{n+1}$.
The discrete surface index is defined by $\bar{j}^n:=\bar{j}(t_{n})$ and we let $\bar{z}^n:=\bar{z}(t_n)$.
For the volumetric flows, we define the averages 
\begin{align*}
Q_\mathrm{f}^n:=\frac{1}{ \tau }\int_{t_n}^{t_{n+1}}Q_\mathrm{f}(t)\,\rmd t
\end{align*}
and assume for simplicity that any of the volumetric flows changes sign at most at the discrete time points $t_n$.
This implies that $\bar{z}(t)$ is monotone on every interval $[t_n,t_{n+1}]$.
To ensure that the surface does not travel more than one cell width~$h$ during $ \tau $, the CFL condition has to imply (cf.\ \eqref{eq:zbarprime})
\begin{equation} \label{eq:zbarprimemax}
 \tau  \max_{0\leq t\leq T}|\bar{z}'(t)|\leq
 \tau \max_{0\leq t\leq T, \atop 0\leq z\leq B}\left\{\frac{|\Qu(t)-\Qf(t)|}{A(z)}, \frac{|\Qu(t)+\Qe(t)|}{A(z)}\right\}
\leq h.
\end{equation} 
To show that the cell concentrations~$\smash{\bC_{\bar{j}^n}^n}$ do not exceed the maximal one~$\hat{X}$, we also introduce the concentration~$\smash{\boldsymbol{\bar{C}}_{\bar{j}^n}^n}$ obtained when all the mass in the surface cell~${\bar{j}^n}$ is located below the surface within the cell; cf.\ Figure~\ref{fig:fig_cells}~(a).
The mass in the cell is
\begin{equation} \label{eq:Cbar}
\boldsymbol{\bar{C}}_{\bar{j}^n}^n A_{\bar{j}^n}\alpha^{n}h
= \bC_{\bar{j}^n}^n A_{\bar{j}^n}h,
\quad\text{where}\quad
\alpha^{n}h := z_{\bar{j}^n+1/2}-\bar{z}^n.
\end{equation} 
We set $\bar{X}_j^n:=\bar{C}_j^{(1),n}+\ldots +\bar{C}_j^{(k_{\bC}),n}$.
Integrating~\eqref{eq:Vprime} from~$t_n$ to $t_{n+1}$, one obtains
\begin{equation} \label{eq:Vdiff}
V \bigl(\bar{z}(t_{n+1}) \bigr) - V \bigl(\bar{z}(t_{n}) \bigr)
= (\bar{Q}^n - \Qu^n) \tau .
\end{equation} 
If the surface stays within one cell between $t_n$ and $t_{n+1}$;  then~\eqref{eq:Vdiff} is equivalent to
\begin{equation} \label{eq:alphadiff}
A_j(\alpha^{n+1} - \alpha^{n})h
= (\bar{Q}^n - \Qu^n) \tau 
\end{equation} 
(cf.\  Figure~\ref{fig:fig_cells}~(b) and (e)). 
During extraction, the surface cannot rise and is thus located somewhere in  cells $\bar{j}^n$ and $\bar{j}^{n}+1$ (cf.\  Figure~\ref{fig:fig_cells}~(d) and (e)).
In light of~\eqref{eq:PhiC_e} and \eqref{eq:PhiS_e}, we approximate the fluxes, which have  to be non-positive, just below the surface in the following way:
\begin{align}
\bPhi^{\bC,n}_{\mathrm{e},\bar{j}^n+1/2}
&:= \biggl(A_{\bar{j}^n+1/2}
\biggl(\vhs(X_{\bar{j}^{n}+1}^n) - \frac{D(X_{\bar{j}^{n}+1}^n)}{h}
\biggr)-{\Qe^n}\biggr)^-
\bC_{\bar{j}^{n}+1}^n,
\label{eq:numflux_eC}\\
%%%%%%%%%%%%%
\bPhi^{\bS,n}_{\mathrm{e},\bar{j}^n+1/2}
&:=  \biggl(
- 
\frac{A_{\bar{j}^n+1/2}X_{\bar{j}^{n}+1}^n}
{\rho_X-X_{\bar{j}^{n}+1}^n}
\biggl(
\vhs(X_{\bar{j}^{n}+1}^n) - \frac{D(X_{\bar{j}^{n}+1}^n)}{h}
\biggr)
- \Qe^n
\biggr)^-
\bS_{\bar{j}^{n}+1}^n.
\label{eq:numflux_eS}
\end{align}

\subsection{Derivation of update formulas}\label{sec:updates}

We here derive the update formulas for~$\smash{\bC_{j}^n}$.
Analogous formulas hold for $\smash{\bS_{j}^n}$ when replacing~$\bC$ by~$\bS$; however, with different definitions of velocities and fluxes.
First come cells that lie below the surface $z=\bar{z}(t)$ at~$t_n$ and~$t_{n+1}$.
Special treatment is needed for the cells near the surface.
All cells strictly above the surface have zero concentrations.
Let $\kappa:= \tau /h$.

\subsubsection*{Cells away from the surface}

Using the integrated form of the balance law on a rectangle $[z_{j-1/2},z_{j+1/2}]\times[t_n,t_{n+1}]$ strictly below the surface; see the gray rectangle in Figure~\ref{fig:fig_cells}~(a), we get the update formula (mass per $h$)
\begin{equation}\label{eq:cellbelow}
A_j\bC_j^{n+1}  = A_j\bC_j^n +  \kappa\bigl(-[\Delta \bPhi^{\bC}]^n_j 
+ hA_j\bR_{\bC,j}^n\bigr),
\end{equation}
and the analogous one for $\smash{\bS_{j}^n}$.
For cell~$N+1$, we get in the similar way the update formula for the underflow concentration; see Section~\ref{subsec:scheme}.

\subsubsection*{Cells near the surface during fill ($t\in T_\mathrm{f}$)}

To obtain a monotone scheme with an invariant-region property, we determine how the mass in the surface cell and the one below evolves; see the trapezoids in Figure~\ref{fig:fig_cells}.
The mass per $h$ at $t_n$ is (we use~\eqref{eq:Cbar})
\begin{equation} \label{eq:initialmass}
\boldsymbol{m}_{\bar{j}^n}^{\bC,n}
= A_{\bar{j}^{n}}\alpha^{n}\boldsymbol{\bar{C}}_{\bar{j}^n}^n 
+ A_{\bar{j}^{n}+1}\bC_{\bar{j}^{n}+1}^n
= A_{\bar{j}^{n}}\bC_{\bar{j}^{n}}^n
+ A_{\bar{j}^{n}+1}\bC_{\bar{j}^{n}+1}^n.
\end{equation} 
During $ \tau $, the feed source along the moving surface is $\Qf^n\bCf^n$, whereas the outflux is 
 $\smash{\bPhi^{\bC,n}_{\bar{j}^n+3/2}}$.
 Thus, by the balance law on any trapezoid  the mass (per~$h$) at~$t_{n+1}$ is
\begin{equation} \label{eq:finalmass}
\boldsymbol{m}_{\bar{j}^n}^{\bC,n+1}
= \boldsymbol{m}_{\bar{j}^n}^{\bC,n}
+ \kappa\bPsi_{\mathrm{f},\bar{j}^{n}}^{\bC,n}
\end{equation} 
where the in- and outflux and source terms are 
\begin{equation*} 
\bPsi_{\mathrm{f},\bar{j}^{n}}^{\bC,n}
:=\Qf^n\bCf^n - \bPhi^{\bC,n}_{\bar{j}^n+3/2} 
+ h\biggl(
\alpha^n A_{\bar{j}^{n}} \bR_{\bC} \biggl( \frac{\bC^{n}_{\bar{j}^n}}{\alpha^n}, \frac{\bS^{n}_{\bar{j}^n}}{\alpha^n} \biggr)
+ A_{\bar{j}^{n}+1}\bR_{\bC,\bar{j}^n+1}^n
\biggr)
\end{equation*} 
(the concentration in cell~$\bar{j}^n$ below the surface is $\smash{\boldsymbol{\bar{C}}_{\bar{j}^n}^n}$). 

\emph{Case (a): Fill case $\bar{j}^{n}=\bar{j}^{n+1}-1$, Figure~\ref{fig:fig_cells}~(a)}:
The surface moves downwards and crosses a cell boundary.
All the mass ends up in one cell: $\smash{A_{\bar{j}^{n}+1}\bC_{\bar{j}^{n}+1}^{n+1} =\boldsymbol{m}_{\bar{j}^n}^{\bC,n+1}}$.

\emph{Case (b): Fill case $\bar{j}^{n}=\bar{j}^{n+1}$, Figure~\ref{fig:fig_cells}~(b)}:
When the surface does not cross any cell boundary during $ \tau $, the mass~\eqref{eq:finalmass} is distributed among the two cells with respect to their volumes (below the surface):
\begin{equation} \label{eq:distrib}
\begin{aligned}
A_{\bar{j}^{n}}\bC_{\bar{j}^{n}}^{n+1} &= \frac{\alpha^{n+1}A_{\bar{j}^{n}}}{
\alpha^{n+1}A_{\bar{j}^{n}}
+ A_{\bar{j}^{n}+1} }
\boldsymbol{m}_{\bar{j}^n}^{\bC,n+1}, 
\\
A_{\bar{j}^{n}+1}\bC_{\bar{j}^{n}+1}^{n+1} &= \frac{A_{\bar{j}^{n}+1}}{
\alpha^{n+1}A_{\bar{j}^{n}}
+ A_{\bar{j}^{n}+1} } \boldsymbol{m}_{\bar{j}^n}^{\bC,n+1}.
\end{aligned}
\end{equation} 

\emph{Case (c): Fill case $\bar{j}^{n}=\bar{j}^{n+1}+1$, Figure~\ref{fig:fig_cells}~(c)}:
After the balance law is used on the purple trapezoid, the final mass (per $h$) is distributed among three cells:
\begin{align*}
A_{\bar{j}^{n}-1}\bC_{\bar{j}^{n}-1}^{n+1} &= \frac{\alpha^{n+1}A_{\bar{j}^{n}-1}}{
\alpha^{n+1}A_{\bar{j}^{n}-1}
+ A_{\bar{j}^{n}} 
+ A_{\bar{j}^{n}+1}}
\boldsymbol{m}_{\bar{j}^n}^{\bC,n+1}, 
\\
A_{\bar{j}^{n}}\bC_{\bar{j}^{n}}^{n+1} &= \frac{A_{\bar{j}^{n}}}{\alpha^{n+1}A_{\bar{j}^{n}-1}
+ A_{\bar{j}^{n}} 
+ A_{\bar{j}^{n}+1}}
\boldsymbol{m}_{\bar{j}^n}^{\bC,n+1}, 
\\
A_{\bar{j}^{n}+1}\bC_{\bar{j}^{n}+1}^{n+1} &= \frac{A_{\bar{j}^{n}+1}}{\alpha^{n+1}A_{\bar{j}^{n}-1}
+ A_{\bar{j}^{n}} 
+ A_{\bar{j}^{n}+1}} \boldsymbol{m}_{\bar{j}^n}^{\bC,n+1}.
\end{align*}

\subsubsection*{Cells near the surface during extraction ($t\in T_\mathrm{e}$)}

During extraction, the surface necessarily moves downwards.
The initial mass is~\eqref{eq:initialmass} and the balance law on a red trapezoid (Figures~\ref{fig:fig_cells}~(d) and~(e)) gives
$ \smash{\boldsymbol{m}_{\bar{j}^n}^{\bC,n+1}
= \boldsymbol{m}_{\bar{j}^n}^{\bC,n}
+ \kappa\bPsi_{\mathrm{e},\bar{j}^n}^{\bC,n},}$
where
\begin{equation*} 
\bPsi_{\mathrm{e},\bar{j}^n}^{\bC,n}
:= \bPhi_{\mathrm{e},\bar{j}^{n}+1/2}^{\bC,n} - \bPhi^{\bC,n}_{\bar{j}^n+3/2}
+ h\biggl(
\alpha^nA_{\bar{j}^n}\bR_{\bC} \biggl( \frac{\bC_{\bar{j}^n}^n}{\alpha^n}, \frac{\bS_{\bar{j}^n}^n}{\alpha^n} \biggr)
+ A_{\bar{j}^n+1}\bR_{\bC,\bar{j}^n+1}^n
\biggr).
\end{equation*} 

\emph{Extraction case $\bar{j}^{n}=\bar{j}^{n+1}-1$, Figure~\ref{fig:fig_cells}~(d)}:
All the mass ends up in one cell: $\smash{A_{\bar{j}^{n}+1}\bC_{\bar{j}^{n}+1}^{n+1} = 
\boldsymbol{m}_{\bar{j}^n}^{\bC,n+1}}$.

\emph{Extraction case $\bar{j}^{n}=\bar{j}^{n+1}$, Figure~\ref{fig:fig_cells}~(e)}:
The surface stays in one cell and we distribute the mass~$\smash{\boldsymbol{m}_{\bar{j}^n}^{\bC,n+1}}$ into two cells with~\eqref{eq:distrib}.

\subsubsection*{The cell in the extraction pipe}
The conservation law for the cell on the $x$-axis gives the mass equality (Figure~\ref{fig:fig_cells}~(d))
\begin{equation*}
A_\mathrm{e}\Delta x\bCe^{n+1} 
= A_\mathrm{e}\Delta x\bCe^{n} 
+  \tau \big(
- \Qe^n\bCe^n - \bPhi^{\bC,n}_{\mathrm{e},\bar{j}^n+1/2}
\big),
\end{equation*}
where the cross-sectional area~$A_\mathrm{e}$ of the effluent pipe is of less importance, since we are only interested in $\bCe^n$ and therefore may choose any $\Delta x$; we set $A_\mathrm{e}\Delta x:= A_1h$.

\subsection{Explicit fully discrete scheme} \label{subsec:scheme}

Given data at $t_n$ and the values $\bar{z}^{n+1}$ and $\bar{j}^{n+1}$, the update formulas for the particulate concentrations are given here and we distinguish between fill and extraction.
We define 
$\lambda_{j}:=\kappa/A_j= \tau /(A_jh)$ and
\begin{alignat*}2
{T}_\mathrm{e}^\mathrm{num} &:=\big\{t_n:\Qe^n>0,\Qf^n=0\big\},
&\qquad
\eta^{n+1}
&:= \frac{A_{\bar{j}^{n}}}{
\alpha^{n+1}A_{\bar{j}^{n}}
+ A_{\bar{j}^{n}+1} },
\\%%%%%%%%
{T}_\mathrm{f}^\mathrm{num} &:=\big\{t_n:\Qe^n=0,\Qf^n\geq0\big\},
&\qquad
\theta^{n+1}
&:= \frac{A_{\bar{j}^{n}}}{
\alpha^{n+1}A_{\bar{j}^{n}-1}
+ A_{\bar{j}^{n}}
+ A_{\bar{j}^{n}+1} }.
\end{alignat*}
The update formulas for the top cells below the surface are different depending on whether fill or extraction occurs.

\subsubsection*{Update formulas for top cells during fill}

If $t_n\in {T}_\mathrm{f}^\mathrm{num}$, then $\bCe^{n+1}=\bzero$. 
The numerical flux is given by~\eqref{eq:numflux_eC} and the marching formulas depending on each layer are computed as follows for the top cells $j \in \{ \bar{j}^n-1,\bar{j}^n, \bar{j}^n+1\}$ (analogously for~$\bS$):
\begin{align*}
 \bC_{j}^{n+1}  & = \omega^{\mathrm{f},n}_j \biggl(\bC_{\bar{j}^{n}}^n
+  \dfrac{A_{\bar{j}^{n}+1} }{A_{\bar{j}^{n}}} \bC_{\bar{j}^{n}+1}^n
+ \lambda_{\bar{j}^n}\bigl(\Qf^n\bCf^n - \bPhi^{\bC,n}_{\bar{j}^n+3/2}\bigr) \\
& \qquad \quad+ 
\tau\biggl(\alpha^n\bR_{\bC} \biggl( \frac{\bC_{\bar{j}^n}^n}{\alpha^n}, \frac{\bS_{\bar{j}^n}^n}{\alpha^n} \biggr) + \dfrac{A_{\bar{j}^n+1}}{A_{\bar{j}^n}}\bR^n_{\bC,\bar{j}^n+1} \biggr)
\biggr)\\
& = \omega^{\mathrm{f},n}_j \biggl(\bC_{\bar{j}^{n}}^n
+  \dfrac{A_{\bar{j}^{n}+1} }{A_{\bar{j}^{n}}} \bC_{\bar{j}^{n}+1}^n
+ \lambda_{\bar{j}^n}\bPsi_{\mathrm{f},\bar{j}^{n}}^{\bC,n}\biggr),
\end{align*}
where the coefficients $\smash{\omega^{\mathrm{f},n}_j}$ are defined in the cases described in Figure~\ref{fig:fig_cells} by:
\begin{alignat*}{5} 
&\mbox{Case (a): $\bar{j}^n = \bar{j}^{n+1}-1$}: 
&&\quad \omega_{\bar{j}^n-1}^{\mathrm{f},n} := 0, 
&&\quad \omega_{\bar{j}^n}^{\mathrm{f},n}   := 0, 
&&\quad \omega_{\bar{j}^n+1}^{\mathrm{f},n} :=  A_{\bar{j}^{n}} / A_{\bar{j}^{n}+1},\\
&\mbox{Case (b): $\bar{j}^n = \bar{j}^{n+1}$}: 
&&\quad \omega_{\bar{j}^n-1}^{\mathrm{f},n} :=0, 
&&\quad \omega_{\bar{j}^n}^{\mathrm{f},n} := \alpha^{n+1}\eta^{n+1}, 
&&\quad \omega_{\bar{j}^n+1}^{\mathrm{f},n} := \eta^{n+1},\\
&\mbox{Case (c): $\bar{j}^n = \bar{j}^{n+1}+1$}: 
&&\quad \omega_{\bar{j}^n-1}^{\mathrm{f},n} := \alpha^{n+1}\theta^{n+1},
&&\quad \omega_{\bar{j}^n}^{\mathrm{f},n}   := \theta^{n+1},
&&\quad \omega_{\bar{j}^n+1}^{\mathrm{f},n} := \theta^{n+1}. 
\end{alignat*}

\subsubsection*{Update formulas for top cells during extraction}
If $t_n\in {T}_\mathrm{e}^\mathrm{num}$, then we compute the numerical fluxes with~\eqref{eq:numflux_eS} and
for the top layers $j \in \{ \bar{j}^n, \bar{j}^n+1\}$, the formula is (analogously for~$\bS$)
\begin{align*}
 \bC_{j}^{n+1}  & = \omega^{\mathrm{e},n}_j \biggl(\bC_{\bar{j}^{n}}^n
+  \dfrac{A_{\bar{j}^{n}+1} }{A_{\bar{j}^{n}}} \bC_{\bar{j}^{n}+1}^n
+ \lambda_{\bar{j}^n}\bigl(\bPhi_{\mathrm{e},\bar{j}^{n}+1/2}^{\bC,n} - \bPhi^{\bC,n}_{\bar{j}^n+3/2}\bigr) \\
& \qquad \quad+ 
\tau\biggl(\alpha^n\bR_{\bC} \biggl( \frac{\bC_{\bar{j}^n}^n}{\alpha^n}, \frac{\bS_{\bar{j}^n}^n}{\alpha^n} \biggr) + \dfrac{A_{\bar{j}^n+1}}{A_{\bar{j}^n}}\bR^n_{\bC,\bar{j}^n+1} \biggr)
\biggr)\\
& =  \omega^{\mathrm{e},n}_j \biggl(\bC_{\bar{j}^{n}}^n
+  \dfrac{A_{\bar{j}^{n}+1} }{A_{\bar{j}^{n}}} \bC_{\bar{j}^{n}+1}^n
+ \lambda_{\bar{j}^n}\bPsi_{\mathrm{e},\bar{j}^{n}}^{\bC,n}\biggr), 
\end{align*}
where the coefficients $\smash{\omega^{\mathrm{e},n}_j}$ are defined in the cases described in Figure~\ref{fig:fig_cells} by:
\begin{alignat*}{4} 
&\mbox{Case (d): $\bar{j}^n=\bar{j}^{n+1}$}:  
&&\quad \omega^{\mathrm{e},n}_{\bar{j}^n}  := \alpha^{n+1}\eta^{n+1}, 
&&\quad \omega^{\mathrm{e},n}_{\bar{j}^n+1}:= \eta^{n+1},\\
&\mbox{Case (e): $\bar{j}^n=\bar{j}^{n+1}-1$}:  
&&\quad \omega^{\mathrm{e},n}_{\bar{j}^n}:= 0,
&&\quad \omega^{\mathrm{e},n}_{\bar{j}^n+1}:= {A_{\bar{j}^{n}}}/{A_{\bar{j}^{n}+1}}.
\end{alignat*}
The effluent concentration is given by
\begin{align*}
 \bCe^{n+1} & = 
\bCe^{n} - \lambda_1
\big( \Qe^n\bCe^{n} 
+ \bPhi^{\bC,n}_{\mathrm{e},\bar{j}^n+1/2}
\big).
\end{align*}

\subsubsection*{Other concentrations}
For the cells $j = \bar{j}^n+2,\dots,N$, the update formula is (analogously for~$\bS$), at every time point $t_n$,
\begin{align}
 \bC_{j}^{n+1} &= \bC_j^n - \lambda_j [\Delta \bPhi^{\bC}]^n_j + \tau \bR_{\bC,j}^n,\label{eq:Cupdate}\\
\bC_{\mathrm{u}}^{n+1} & = \bC_{\mathrm{u}}^n + \lambda_{N+1} \bigl(\bPhi^{\bC,n}_{N+1/2} - \Qu^n\bC_{\mathrm{u}}^n\bigr).\label{eq:Cunum}
\end{align}
where the numerical flux is computed by~\eqref{eq:numfluxes}.
For the cells above the surface, we have $\bC_{j}^{n+1} = \boldsymbol{0}$ for $j<\bar{j}^n-1$.
Finally, one computes
\begin{align*}
X_j^{n} = C_j^{(1),n}+\cdots+C_j^{(k_{\bC}),n},
\qquad
W_j^{n} = \rho_L(1-X_j^n/\rho_X) - \bigl(S_j^{(1),n}+\cdots + S_j^{(k_{\bS}),n}\bigr).
\end{align*}

\subsection{A splitting scheme and invariant-region property} \label{subsec:numcfl}

It is desirable that the solution vectors $\bU:=(\bC,\bS)$ and $\bU_\mathrm{e}:=(\bCe,\bSe)$ of the model~\eqref{finalmod} stay in the set
\begin{equation*} 
\Omega := \bigl\{ \bU \in\mathbb{R}^{k_{\bC}+k_{\bS}}: \bC\geq\bzero, \, \bS\geq\bzero, \, 
C^{(1)} + \dots + C^{(k_{\bC})}\leq\hat{X}
\bigr\}.
\end{equation*}
To ensure that~$\Omega$ is an invariant set for the numerical solutions, we split the scheme in Section~\ref{subsec:scheme} by taking one time step without the reaction terms and then one time step with only the reaction terms.
For a cell strictly below the surface and the solid concentrations, the update formula~\eqref{eq:Cupdate} can be written
\begin{equation*}
\bC_j^{n+1} = \bigl\{\bC_j^n - \lambda_j[\Delta \bPhi^{\bC}]^n_j\bigr\}
+ \tau\bR_{\bC}(\bC^{n}_j,\bS^{n}_j).
\end{equation*}
and similarly for $\bS_j^n$.
The splitting principle is to compute the expression in the curled brackets first and then use the result for the second step as follows:
\begin{alignat}2
\boldsymbol{\check{\bC}}_j^{n+1} &= \bC_j^n - \lambda_j[\Delta \bPhi^{\bC}]^n_j,&\qquad&
\boldsymbol{\check{\bS}}_j^{n+1} = \bS_j^n - \lambda_j[\Delta \bPhi^{\bS}]^n_j,\label{eq:Ccheck}\\
\bC_j^{n+1} &= \boldsymbol{\check{\bC}}_j^{n+1}
+ \tau\bR_{\bC}(\boldsymbol{\check{\bC}}_j^{n+1}, \boldsymbol{\check{\bS}}_j^{n+1}),
&\qquad&
\bS_j^{n+1} = \boldsymbol{\check{\bS}}_j^{n+1}
+ \tau\bR_{\bS}(\boldsymbol{\check{\bC}}_j^{n+1}, \boldsymbol{\check{\bS}}_j^{n+1}).
\label{eq:ODEstep}
\end{alignat}
For the cells involving the surface, the first-step update formulas, $j\in\{\bar{j}^n-1,\bar{j}^n,\bar{j}^n+1\}$ are
\begin{alignat}2
 \check{\bC}_{j}^{n+1}  & = \omega^{\mathrm{f},n}_j \biggl(\bC_{\bar{j}^{n}}^n
+  \frac{A_{\bar{j}^{n}+1} }{ A_{\bar{j}^{n}}}  \bC_{\bar{j}^{n}+1}^n
+ \lambda_{\bar{j}^n}\bigl(\Qf^n\bCf^n - \bPhi^{\bC,n}_{\bar{j}^n+3/2}
\bigr) \biggr) &\quad& \text{if $t_n\in T_\mathrm{f}^\mathrm{num}$,} \label{eq:Ccheck_f}
\\
 \check{\bC}_{j}^{n+1}  & = \omega^{\mathrm{e},n}_j \biggl(\bC_{\bar{j}^{n}}^n
+  \frac{A_{\bar{j}^{n}+1} }{ A_{\bar{j}^{n}}}  \bC_{\bar{j}^{n}+1}^n
+ \lambda_{\bar{j}^n}\bigl(\bPhi_{\mathrm{e},\bar{j}^{n}+1/2}^{\bC,n} - \bPhi^{\bC,n}_{\bar{j}^n+3/2}
\bigr) \biggr) 
&\quad& \text{if $t_n\in T_\mathrm{e}^\mathrm{num}$,} \label{eq:Ccheck_e}
\end{alignat}
and the analogous formulas with $\bC$ replaced by~$\bS$.
The second step consists in the formulas
\begin{alignat}2
 \bC_{j}^{n+1}  & = \check{\bC}_{j}^{n+1} + \omega^{\mathrm{f},n}_{j} \tau\biggl(\alpha^n \bR_{\bC}\biggl(\frac{\check{\bC}_{\bar{j}^n}^{n+1}}{\alpha^n}, \frac{\check{\bS}_{\bar{j}^n}^{n+1}}{\alpha^n}\biggr) + \dfrac{A_{\bar{j}^n+1}}{A_{\bar{j}^n}} \bR_{\bC}\bigl(\check{\bC}_{\bar{j}^n+1}^{n+1}, \check{\bS}_{\bar{j}^n+1}^{n+1}\bigr) \biggr) 
&\quad& \text{if $t_n\in T_\mathrm{f}^\mathrm{num}$,} \label{eq:ODEstep_f}\\
 \bC_{j}^{n+1}  & = \check{\bC}_{j}^{n+1} + \omega^{\mathrm{e},n}_{j} \tau\biggl(\alpha^n \bR_{\bC}\biggl(\frac{\check{\bC}_{\bar{j}^n}^{n+1}}{\alpha^n}, \frac{\check{\bS}_{\bar{j}^n}^{n+1}}{\alpha^n}\biggr) + \dfrac{A_{\bar{j}^n+1}}{A_{\bar{j}^n}} \bR_{\bC}\bigl(\check{\bC}_{\bar{j}^n+1}^{n+1}, \check{\bS}_{\bar{j}^n+1}^{n+1}\bigr) \biggr)
&\quad& \text{if $t_n\in T_\mathrm{e}^\mathrm{num}$,} \label{eq:ODEstep_e}
\end{alignat}
and the analogous formulas with $\bC$ replaced by~$\bS$.

The time step~$\tau$ has to be bounded by the  CFL condition 
\begin{equation}\label{cfl}\tag{CFL}
 \tau \, \max\big\{\beta_1,\beta_2,M_{\bC}(1+M_3),M_{\bS},\tilde{M}/\varepsilon\big\}
\leq 1,
\end{equation} 
where $\beta_1$ and $\beta_2$ depend on~$h$, $h^2$, the volumetric flows, and the constitutive functions by
\begin{align*}
\beta_1 &:= \frac{\|Q\|_{T}}{A_\mathrm{min}h}  + \frac{M_1}{h}
\big(
\|\vhs'\|\hat{X} + \vhs(0)
\big) 
+ \frac{2M_2}{h^2}
\big(
\|d\|\hat{X} + D(\hat{X})
\big)
, \\
\beta_2 &:= 
\max\{M_1,1\}
\dfrac{\rho_X+\hat{X}}{\rho_X-\hat{X}} \frac{\|Q\|_{T}}{A_\mathrm{min}h} 
+ \dfrac{\hat{X} M_1}{\rho_X-\hat{X}}\dfrac{2\vhs(0)}{h}
+ \dfrac{\hat{X} M_2}{\rho_X-\hat{X}}\dfrac{D(\hat{X})}{h^2},
\end{align*}
and where the constants are given by (here, $\xi$ represents $\vhs, \vhs'$ or $d$) (cf.~\eqref{eq:positivity})
\begin{align*}
&M_{\boldsymbol{\xi}}:=\sup_{\boldsymbol{\mathcal{U}}\in\Omega, \atop 1\le k\le k_{\boldsymbol{\xi}}}
\sum_{l\in I_{\boldsymbol{\xi},k}^-}|\sigma_{\boldsymbol{\xi}}^{(k,l)}|\bar{r}_{\boldsymbol{\xi}}^{(l)} (\bC,\bS),
\quad\boldsymbol{\xi}\in\{\bC,\bS\},
\qquad 
\tilde{M}:=\sup_{\boldsymbol{\mathcal{U}}\in\Omega}
{\tilde{R}_{\bC}}(\bC,\bS),
 \\
%%%
&\|\xi\|:=\max\limits_{0\le X\le\hat{X}}|\xi(X)|, \qquad\|Q\|_{T}:=\max_{0\le t\le T}{\big\{|\Qu(t)-\Qf(t)|,\Qu(t)+\Qe(t)\big\}},\\
& M_1:= \underset{j=1,\dots,N}{\max}\left\{\dfrac{A_{j+1/2}}{A_j},\dfrac{A_{j-1/2}}{A_j}\right\},\qquad
    M_2:= \underset{j=1,\dots,N}{\max}\left\{\dfrac{A_{j+1/2}+A_{j-1/2}}{A_j}\right\},\\
&M_3:=\underset{j=1,\dots,N-1}{\max}\left\{\dfrac{A_{j}}{A_{j+1}}\right\}.
\end{align*}

\begin{theorem}\label{theorem}
Consider the numerical splitting method in Section~\ref{subsec:numcfl}.
If $\bU_j^n:=(\bC_j^n,\bS_j^n)\in\Omega$ for all $j\ne\bar{j}^n$, $\boldsymbol{\bar{\mathcal{U}}}_{\bar{j}^n}^n:=(\boldsymbol{\bar{C}}_{\bar{j}^n}^n,\boldsymbol{\bar{S}}_{\bar{j}^n}^n)\in\Omega$, $\bU_{\mathrm{e}}^n:=(\bCe^n,\bSe^n)\in\Omega$ and \eqref{cfl} holds, then 
\begin{align*} 
\bU_j^{n+1},\boldsymbol{\bar{\mathcal{U}}}_{\bar{j}^n}^{n+1}, \bU_{\mathrm{e}}^{n+1}\in\Omega.
\end{align*} 
\end{theorem}
Theorem~\ref{theorem} is proved by the following lemmas and by the fact that~\eqref{cfl} implies~\eqref{eq:zbarprimemax}, which we have used in the derivation of the scheme.
For the proofs write the update formulas as 
\begin{align*} 
C^{(k),n+1}_j= \calH_{\bC,j}^{(k),n} \bigl(\bC_{j-1}^{n},\dots,\bC_{j+3}^{n},\bS_j^n \bigr)
\end{align*} 
for one component~$k\in\{1,\ldots,k_{\bC}\}$ (see Section~\ref{subsec:scheme}).
That formula includes the underflow concentrations~\eqref{eq:Cunum} for $j=N+1$.  
Summing for fixed~$j$  all components of the update formula~\eqref{eq:Cupdate} for $\smash{\bC_j^{n+1}}$, one gets
\begin{equation}\label{eq:Xupdate}
X_j^{n+1}=X_j^{n}+\lambda_j[\Delta \tilde{\Phi}^{\bC}]^n_j + \tau\tilde{R}_{\bC,j}^n,
\end{equation}
where $\smash{\tilde{\Phi}}^{\bC}$ and $\smash{\tilde{R}_{\bC,j}^n}$ denote the sum of all components, respectively.
Analogous considerations lead to formulas for~$\bar{X}_j^n$ and $X_\mathrm{e}^{n+1}$.

\begin{lemma}\label{lem:X}
Let $\bR_{\bC}\equiv\bzero$.
If $\bU_j^n:=(\bC_j^n,\bS_j^n)\in\Omega$ for all $j\ne\bar{j}^n$, $\boldsymbol{\bar{\mathcal{U}}}_{\bar{j}^n}^n:=(\boldsymbol{\bar{C}}_{\bar{j}^n}^n,\boldsymbol{\bar{S}}_{\bar{j}^n}^n)\in\Omega$, $\bU_{\mathrm{e}}^n:=(\bCe^n,\bSe^n)\in\Omega$ and \eqref{cfl} holds, then 
\begin{align*} 
0\leq X_j^{n+1},{\bar{X}}_{\bar{j}^n}^{n+1},X_\mathrm{e}^{n+1} \leq\hat{X} \quad \text{\em for all~$j$}.
\end{align*} 
\end{lemma}

\begin{proof}
We write the general update formula~\eqref{eq:Xupdate} as $X_j^{n+1}=\calH_{X,j}^{n}(X_{j-1}^n,\ldots,X_{j+3}^n)$ and let this include the surface concentration~$\bar{X}_{\bar{j}^n}^n$.
For cells away from the moving surface, we refer to~\cite[Theorem~3.1]{bcdimamat21} from which we also collect
\begin{equation*}
v_{j+1/2}^{X,n,+}
= 
\bigl(
q_{j+1/2}^n + \gamma_{j+1/2}\big(
\vhs(X_{j+1}^n) - J_{j+1/2}^{\bC,n}\big)
\bigr)^+
\leq \frac{\Qu^n}{A_{j+1/2}} 
+ \vhs(0) + \frac{D(\hat{X})}{h}.
\end{equation*}
To show the monotonicity for the cells near the surface, we set 
\begin{align*}
\chi^+&:=\chi_{\{{v}^{X,n}_{\bar{j}^{n}+1/2}\geq 0\}},
\qquad
\chi^-:=\chi_{\{{v}^{X,n}_{\bar{j}^{n}+1/2}\leq 0\}},
\qquad
\chi_{\mathrm{e}}^-:=\chi_{\{v_{\mathrm{e},\bar{j}^{n}+1/2}^{X,n}\leq 0\}},
\\
\nu_1&:= \hat{X}\left(
 \frac{\|d\|}{h} 
+ \vhs(0)
\right)
+  \frac{D(\hat{X})}{h},
\qquad
\nu_2:= 
\hat{X} \left(\frac{\|d\|}{h} + \|\vhs'\|\right) 
+  \frac{D(\hat{X})}{h},
\end{align*}
and calculate for $j=\bar{j}^n$
\begin{align*}
\pp{{v}^{X,n,\pm}_{j+1/2}}{X_j^n}
&= \chi^{\pm}\frac{d({X}_{j}^n)}{h},
\qquad
\pp{{v}^{X,n,\pm}_{j+1/2}}{X_{j+1}^n}
= \vhs'(X_{j+1}^n) 
- \chi^{\pm}\frac{d(X_{j+1}^n)}{h},
\\
%%%%%%%%
\pp{{\tilde{\Phi}}^{\bC,n}_{j+1/2}}{X_j^n}
& = A_{j+1/2}\biggl(
\left(
\chi^-X_{j+1}^n 
+ \chi^+{X}_{j}^n
\right)
\frac{d({X}_j^n)}{h} 
+ v^{X,n,+}_{j+1/2}
\biggr)\leq A_{j+1/2}\nu_1
+ \Qu^n,
\\
%%%%%%%
\pp{{\tilde{\Phi}}^{\bC,n}_{j+1/2}}{X_{j+1}^n}
&= A_{j+1/2}\biggl(
{v}^{X,n,-}_{j+1/2}
+ \left(
\chi^-X_{j+1}^n
+ \chi^+{X}_j^n\right)
\left(\vhs'(X_{j+1}^n) - \frac{d(X_{j+1}^n)}{h}
\right)
\biggr)\leq 0,
\\
%%%%%%%%
\pp{\tilde{\Psi}_{\mathrm{f},j}^{\bC,n}}{X_{j-1}^n} 
&= 0,\qquad
\pp{\tilde{\Psi}_{\mathrm{f},j}^{\bC,n}}{X_{j}^n}
= 0,
\\
%%%%%%5
\pp{\tilde{\Psi}_{\mathrm{f},j}^{\bC,n}}{X_{j+1}^n}
&= - \pp{{\tilde{\Phi}}_{j+3/2}^{\bC,n}}{X_{j+1}^n}
\geq -A_{j+3/2}\nu_1 - \Qu^n,
\qquad
\pp{\tilde{\Psi}_{\mathrm{f},j}^{\bC,n}}{X_{j+2}^n}
=- \pp{{\tilde{\Phi}}_{j+3/2}^{\bC,n}}{X_{j+2}^n}
\geq 0.
\end{align*}
For ease of notation, we introduce 
\begin{align}\label{eq:Upsilon_f}
 \boldsymbol{\Upsilon}_{\mathrm{f},j}^{\bC,n} := \bC_{j}^n
+  \dfrac{A_{j+1} }{A_{j}} \bC_{j+1}^n
+ \lambda_{j}\bPsi_{\mathrm{f},j}^{\bC,n}
\end{align}
and let as usual tilde denote the sum of all components of a vector.
In the case $t_n\in {T}_\mathrm{f}^{\rm num}$, all the coefficients~$\omega^{\mathrm{f},n}_j$ for $\smash{\tilde{\Upsilon}_{\mathrm{f},j}^{\bC,n}}$ in~\eqref{eq:Ccheck_f} are non-negative, so it suffices to show that the derivatives of this function are non-negative under~\eqref{cfl}:
\begin{align*}
\pp{\tilde{\Upsilon}_{\mathrm{f},j}^{\bC,n}}
{X_{j-1}^n} 
&=
\lambda_{j}\pp{\tilde{\Psi}_{\mathrm{f},j}^{\bC,n}}{X_{j-1}^n} = 0,
\qquad
\pp{\tilde{\Upsilon}_{\mathrm{f},j}^{\bC,n}}
{X_{j+2}^n} 
= \lambda_{j}\pp{\tilde{\Psi}_{\mathrm{f},j}^{\bC,n}}{X_{j+2}^n}
\geq 0,
\\
%%%%%%%%%%%%%%%
\pp{\tilde{\Upsilon}_{\mathrm{f},j}^{\bC,n}}
{X_{j}^n} 
&= 1 + \lambda_{j}\pp{\tilde{\Psi}_{\mathrm{f},j}^{\bC,n}}{X_{j}^n}
\geq 1 + 0 
\geq 0,\qquad
\pp{\tilde{\Upsilon}_{\mathrm{f},j}^{\bC,n}}
{X_{j+1}^n} 
= \frac{A_{j+1}}{A_{j}} + \lambda_{j}\pp{\tilde{\Psi}_{\mathrm{f},j}^{\bC,n}}{X_{j+1}^n}\\
%%%%%%%%%%%%
&\geq \frac{A_{j+1}}{A_{j}}
- \lambda_{j}\big(
A_{j+3/2}\nu_1 + \Qu^n
\big)
\geq \frac{A_{j+1}}{A_{j}}\left(
1 - { \tau }\left(
\frac{M_1\nu_1}{h} 
+ \frac{\|Q\|_{T}}{A_{\mathrm{min}}h} 
\right)
\right) 
\geq 0.
\end{align*}

In the case $t_n\in {T}_\mathrm{e}^{\rm num}$, we first estimate
\begin{align*}
\pp{\tilde{\Phi}^{\bC,n}_{\mathrm{e},j+1/2}}{X_{j+1}^n}
&= A_{j+1/2}\left(
\chi_\mathrm{e}^-\left(
\vhs'(X_{j+1}^n) - \frac{d(X_{j+1}^n)}{h}
\right)X_{j+1} 
+ v_{\mathrm{e},j+1/2}^{X,n,-}
\right)\geq - A_{j+1/2}\nu_2 -\Qe^n.
\end{align*}
This estimation shall be added to the derivatives of~$\smash{\tilde{\Psi}_{\mathrm{f},j}^{\bC,n}}$ to obtain those for~$\smash{\tilde{\Psi}_{\mathrm{e},j}^{\bC,n}}$.
Introducing
\begin{align}\label{eq:Upsilon_e}
 \boldsymbol{\Upsilon}_{\mathrm{e},j}^{\bC,n} := \bC_{j}^n
+  \dfrac{A_{j+1} }{A_{j}} \bC_{j+1}^n
+ \lambda_{j}\bPsi_{\mathrm{e},j}^{\bC,n},
\end{align}
we see that the only derivative that differs from those above is 
\begin{align*}
\pp{\tilde{\Upsilon}_{\mathrm{e},j}^{\bC,n}}
{X_{j+1}^n} 
&= 
\frac{A_{j+1}}{A_{j}} + \lambda_{j}\pp{\tilde{\Psi}_{\mathrm{e},j}^{\bC,n}}{X_{j+1}^n}
\geq \frac{A_{j+1}}{A_{j}}
- \lambda_{j}\big(
A_{j+1/2}\nu_2 + \Qe^n + A_{j+3/2}\nu_1 + \Qu^n 
\big)
\\
&\geq \frac{A_{j+1}}{A_{j}}\left(
1 - { \tau }\left(
\frac{M_1(\nu_1+\nu_2)}{h} + \frac{\|Q\|_{T}}{A_{\mathrm{min}}h} 
\right)
\right) 
\geq 0.
\end{align*}
We write $C_\mathrm{e}^{(k),n+1} =\calH_{\bCe}^{(k),n}(\bCe^n,\bC_{\bar{j}^n+1}^n)$ for the update formula for one component (they are all equal) of the effluent concentration.
Summing all the components, we obtain
$X_\mathrm{e}^{n+1} =\calH_{\bCe}^{(k),n}(\Xe^n,X_{\bar{j}^n+1}^n)$ for any fixed $k\in\{1,\ldots,k_{\bC}\}$.
This formula is trivial for $t_n\in {T}_\mathrm{f}^{\mathrm{num}}$, and for $t_n\in {T}_\mathrm{e}^{\mathrm{num}}$, we get
\begin{equation*}
\pp{\calH_{\bCe}^{(k),n}}{X_\mathrm{e}^{n}}
= 1 - \lambda_1\Qe^n
\geq 1
-  \tau  \frac{\|Q\|_{T}}{A_{\mathrm{min}}h}
\geq 0,\qquad
\pp{\calH_{\bCe}^{(k),n}}{X_{\bar{j}^n+1}^n}
= - \lambda_1 
\pp{\tilde{\Phi}^{\bC,n}_{\mathrm{e},j+1/2}}{X_{j+1}^n}
\geq 0.
\end{equation*}

To prove the boundedness, the monotonicity in each variable of $\smash{\calH_{X,j}^n}$ and the assumption~\eqref{eq:techrel} are used to obtain, for $t_n\in {T}_\mathrm{f}^{\mathrm{num}}$ and $j=\bar{j}^n=\bar{j}^{n+1}$,
\begin{align*}
0&\leq \alpha^{n+1}\eta_{j}^{n+1} \lambda_j\Qf^n X_{\mathrm{f}}^{n} = \calH_{X,j}^n(0,\dots,0)
\leq X_j^{n+1}
= \calH_{X,j}^n \bigl(
0,X_{j}^{n},\dots,X_{j+3}^{n}
\bigr)\\
&
\leq \calH_{X,j}^n
\bigl(0,\alpha^n{\hat{X}},{\hat{X}},{\hat{X}},{\hat{X}}\bigr)= \alpha^{n+1}\eta_{j}^{n+1}
\biggl(
\alpha^n\hat{X} +  \frac{A_{j+1}}{A_{j}}  \hat{X}
+ \lambda_j \bigl(
\Qf^n X_\mathrm{f}^{n} - \Qu^n\hat{X}
 \bigr)
\biggr)\\
&\leq \alpha^{n+1}\eta_{j}^{n+1}\hat{X}\biggl(
\alpha^n +  \frac{A_{j+1}}{ A_{j} }
+ \lambda_j (
\Qf^n - \Qu^n
)
\biggr) \overset{\eqref{eq:alphadiff}}{=} \alpha^{n+1}\hat{X}\eta_{j}^{n+1}\biggl(
 \frac{A_{j+1}}{A_{j}} 
+ \alpha^{n+1}
\biggr)
= \alpha^{n+1}\hat{X}.
\end{align*}
For $t_n\in {T}_\mathrm{f}^{\mathrm{num}}$ and  $j=\bar{j}^n-1=\bar{j}^{n+1}$, we first see that~\eqref{eq:Vdiff} implies; cf.\ Figure~\ref{fig:fig_cells}(c),
\begin{align} \label{eq:green_alpha}
\alpha^{n+1}A_{\bar{j}^n-1} + A_{\bar{j}^n}
- \alpha^nA_{\bar{j}^n} = (\Qf^n-\Qu^n)\kappa,
\end{align} 
which we use at the end of the following estimate:
\begin{align*}
0&\leq \alpha^{n+1}\theta^{n+1}
\lambda_{\bar{j}^n-1}\Qf X_{\mathrm{f}}^{n} = \calH_{X,\bar{j}^n-1}^n(0,\dots,0)
\leq X_{\bar{j}^n-1}^{n+1}\\
&
= \calH_{X,\bar{j}^n-1}^n\bigl(
0,0,X_{\bar{j}^n}^{n},X_{\bar{j}^n+1}^{n},X_{\bar{j}^n+2}^{n}
\bigr)
\leq \calH_{X,\bar{j}^n-1}^{n}
\bigl(0,0,\alpha^n{\hat{X}},{\hat{X}},{\hat{X}}\bigr)\\
&= \alpha^{n+1}\theta^{n+1}
\biggl(
\alpha^n\hat{X} +  \frac{A_{\bar{j}^n+1}}{ A_{\bar{j}^n}}  \hat{X}
+ \lambda_{\bar{j}^n} \bigl(
\Qf^n X_\mathrm{f}^{n} - \Qu^n\hat{X}
\bigr)
\biggr)\\
&\leq \alpha^{n+1}\hat{X}\theta^{n+1}
\biggl(
\alpha^n +  \frac{A_{\bar{j}^n+1} }{ A_{\bar{j}^n} } 
+ \alpha^{n+1} \frac{A_{\bar{j}^n-1}}{ A_{\bar{j}^n}}  + 1
- \alpha^n
\biggr)
=\alpha^{n+1}\hat{X}.
\end{align*}
In the case $j=\bar{j}^n=\bar{j}^{n+1}+1$ (cf.\ Figure~\ref{fig:fig_cells}(c))  we can still use~\eqref{eq:green_alpha} to obtain
\begin{align*}
0&\leq \theta^{n+1}
\lambda_{\bar{j}^n}\Qf X_{\mathrm{f}}^{n} = \calH_{X,\bar{j}^n}^{n}(0, \dots ,0)
\leq X_{\bar{j}^n}^{n+1}\\
&
= \calH_{X,\bar{j}^n}^{n} \bigl(
0,0,X_{\bar{j}^n}^{n},X_{\bar{j}^n+1}^{n},X_{\bar{j}^n+2}^{n}
\bigr)
\leq \calH_{X,\bar{j}^n}^{n}
\bigl(0,0,\alpha^n{\hat{X}},{\hat{X}},{\hat{X}}\bigr)\\
&= \theta^{n+1}
\biggl(
\alpha^n\hat{X} +  \frac{A_{\bar{j}^n+1}}{ A_{\bar{j}^n}} \hat{X}
+ \lambda_{\bar{j}^n} \bigl(
\Qf^n X_\mathrm{f}^{n} - \Qu^n\hat{X}
\bigr) \biggr)\\
&\leq \hat{X}\theta^{n+1}
\biggl(
\alpha^n +  \frac{A_{\bar{j}^n+1}}{A_{\bar{j}^n}}  
+ \alpha^{n+1} \frac{A_{\bar{j}^n-1}}{ A_{\bar{j}^n}} + 1
- \alpha^n
\biggr)
=\hat{X}.
\end{align*}
A similar estimation can be made for the case $j=\bar{j}^n+1=\bar{j}^{n+1}+2$.
For the case~$j=\bar{j}^n+1=\bar{j}^{n+1}$; see Figure~\ref{fig:fig_cells}(a), we note that~\eqref{eq:Vdiff} implies
\begin{equation*}
\alpha^{n+1}A_{\bar{j}^n+1} - (\alpha^nA_{\bar{j}^n}+A_{\bar{j}^n+1}) = (\Qf^n-\Qu^n)\kappa,
\end{equation*}
which we use to estimate
\begin{equation*}
\begin{aligned}
0&\leq  \frac{A_{\bar{j}^n} }{ A_{\bar{j}^n+1}}  
\lambda_{\bar{j}^n}\Qf X_{\mathrm{f}}^{n}
= \calH_{X,\bar{j}^n+1}^{n}(0,\dots,0)
\leq X_{\bar{j}^n+1}^{n+1}\\
& = \calH_{X,\bar{j}^n+1}^{n} \bigl(0,
X_{\bar{j}^n}^{n},\dots ,X_{\bar{j}^n+3}^{n}
\bigr)\leq \calH_{X,\bar{j}^n+1}^{n}
\bigl(0,
\alpha^n{\hat{X}},{\hat{X}},{\hat{X}},{\hat{X}}
\bigr)\\
&=  \frac{ A_{\bar{j}^n} }{ A_{\bar{j}^n+1}}  \biggl(
\alpha^n\hat{X}
+  \frac{A_{\bar{j}^n+1}}{ A_{\bar{j}^n}}  \hat{X}
+ \lambda_{\bar{j}^n} \bigl(
\Qf^n X_\mathrm{f}^{n} - \Qu^n\hat{X}
 \bigr)
\biggr)\\
&\leq  \frac{A_{\bar{j}^n} }{ A_{\bar{j}^n+1}}  \hat{X}\biggl(
\alpha^n
+  \frac{A_{\bar{j}^n+1}}{ A_{\bar{j}^n} } 
+ \alpha^{n+1} \frac{A_{\bar{j}^n+1} }{A_{\bar{j}^n}}  - \biggl(\alpha^n +  \frac{A_{\bar{j}^n+1}}{A_{\bar{j}^n}} \biggr)
\biggr) =\alpha^{n+1}\hat{X}.
\end{aligned}
\end{equation*}
The remaining fill cases are similar;  we omit details.
For $t_n\in {T}_\mathrm{e}^{\mathrm{num}}$,  similar estimations apply; the only difference is that $\smash{\Qf X_{\mathrm{f}}^{n}}$ is replaced by $\smash{\tilde{\Phi}_{\mathrm{e},j+1/2}^{\bC,n}}$, which equals zero when the concentrations are zero, to prove the lower bound.
For the upper bound, one uses~\eqref{eq:Vdiff} with $\bar{Q}^n=-\Qe^n$ instead of $\Qf^n$.
For the effluent, we get
\begin{align*}
0&
= \calH_{\bC_{\rm e}}^{n}(0,0)
\leq X_{\rm e}^{n+1}
= \calH_{\bC_{\rm e}}^{n} \bigl(
\Xe^{n},X_{\bar{j}^n+1}^{n}
\bigr)
\leq \calH_{\bC_{\rm e}}^{n} \bigl(
{\hat{X}},{\hat{X}}
\bigr) 
\\ & = \hat{X} -\lambda_1(\Qe^n+\Qu^n)\hat{X}
\leq \hat{X}.
\end{align*}
\end{proof}

\begin{lemma}\label{lemmaC}
Let $\bR_{\bC}\equiv\bzero$.
If $\bU_j^n:=(\bC_j^n,\bS_j^n)\in\Omega$ for all $j\ne\bar{j}^n$, $\boldsymbol{\bar{\mathcal{U}}}_{\bar{j}^n}^n:=(\boldsymbol{\bar{C}}_{\bar{j}^n}^n,\boldsymbol{\bar{S}}_{\bar{j}^n}^n)\in\Omega$, $\bU_{\mathrm{e}}^n:=(\bCe^n,\bSe^n)\in\Omega$ and \eqref{cfl} holds, then 
\begin{align*} 
0\leq\bC_j^{n+1},\boldsymbol{\bar{C}}_{\bar{j}^n}^{n+1}, \bCe^{n+1}\leq\hat{X} \quad  \text{\em for all~$j$}.
\end{align*} 
\end{lemma}

\begin{proof}
This follows directly from Lemma~\ref{lem:X}, since each component of the update formula for $\smash{\bC_j^n}$ is equal to that of $\smash{X_j^n}$ (there is no reaction term).
\end{proof}

\begin{lemma}\label{lemmaS}
Let $\bR_{\bS}\equiv\bzero$.
If $\bU_j^n:=(\bC_j^n,\bS_j^n)\in\Omega$ for all $j\ne\bar{j}^n$, $\boldsymbol{\bar{\mathcal{U}}}_{\bar{j}^n}^n:=(\boldsymbol{\bar{C}}_{\bar{j}^n}^n,\boldsymbol{\bar{S}}_{\bar{j}^n}^n)\in\Omega$, $\bU_{\mathrm{e}}^n:=(\bCe^n,\bSe^n)\in\Omega$ and \eqref{cfl} holds, then $\bS_j^{n+1}\geq 0$ for all~$j$ and $\bSe^{n+1}\geq 0$.
\end{lemma}

\begin{proof}
We start as in the proof of Lemma~\ref{lem:X}, use the notation and estimations from there, and prove monotonicity of each component of the right-hand side~$\smash{\calH_{\bS,\bar{j}^n}^{(k),n}}$ of the update formula for component~$\smash{S_j^{(k),n}}$, which we write as~$\smash{S_j^{n}}$.  We also skip the superscript~$(k)$ for components of other vectors.
We let $\hat{\rho}:=1/(\rho_X-\hat{X})$ and
\begin{equation*}
\nu_3:=\hat{\rho}\biggl(
 (\rho_X + \hat{X}
 ) \frac{\|Q\|_{T}}{A_{\mathrm{min}}}  
+ \biggl( \vhs(0) 
+  \frac{D(\hat{X})}{h}\biggr)\hat{X}
\biggr).
\end{equation*}
The numerical fluxes are different and we get
\begin{align*}
\pp{{\Phi}^{\bS,n}_{j+1/2}}{S_j^n}
&= A_{j+1/2}
 \dfrac{(\rho_Xq_{j+1/2}^n-F^{X,n}_{j+1/2})^+}{\rho_X-X_{j}^n}
\leq A_{j+1/2}\hat{\rho}\big(
\rho_X q_{j+1/2}^n
- v_{j+1/2}^{X,n,-}\hat{X}
\big)
\\
&= A_{j+1/2}\hat{\rho}\big(
\rho_X q_{j+1/2}^n
+ v_{j+1/2}^{X,n,+}\hat{X}
\big)\\
&\leq A_{j+1/2}\hat{\rho}\biggl(
\rho_Xq_{j+1/2}^n 
+  \biggl(q_{j+1/2}^n 
+ \vhs(0) 
+  \frac{D(\hat{X})}{h} \biggr)\hat{X}
\biggr)
\leq A_{j+1/2}\nu_3,
\\
%%%%%%%
\pp{{\Phi}^{\bS,n}_{j+1/2}}{S_{j+1}^n}
&= A_{j+1/2}
\dfrac{(\rho_Xq_{j+1/2}^n -F^{X,n}_{j+1/2})^-}{\rho_X-X_{j+1}}
\leq 0.
\end{align*}
Because of the similarities between $\Psi_{\mathrm{f},j}^{\bC,n}$ and $\Psi_{\mathrm{f},j}^{\bS,n}$, we only write the difference here:
\begin{align*}
\pp{\Psi_{\mathrm{f},j}^{\bS,n}}{S_{j+1}^n}
= - \pp{{\Phi}_{j+3/2}^{\bS,n}}{S_{j+1}^n}
\geq -A_{j+3/2}\nu_3.
\end{align*}

To prove the monotonicity in the case $t_n\in {T}_\mathrm{f}^{\rm num}$, we conclude that the estimations are in fact similar to those in the proof of Lemma~\ref{lemmaC} with the following difference (where $\boldsymbol{\Upsilon}_{\mathrm{f},j}^{\bS,n}$ is defined analogously to~\eqref{eq:Upsilon_f} and ${\Upsilon}_{\mathrm{f},j}^{\bS,n}$ denotes an arbitrary component of that vector):
\begin{equation*}
\pp{{\Upsilon}_{\mathrm{f},j}^{\bS,n}}
{S_{j+1}^n} 
\geq \frac{A_{j+1}}{A_{j}}\left(
1 - { \tau }
\frac{M_1\nu_3}{h}  
\right) 
\geq 0.
\end{equation*}

To prove the monotonicity in the case $t_n\in {T}_\mathrm{e}^{\rm num}$, we first estimate
\begin{align*}
\pp{\Phi^{\bS,n}_{\mathrm{e},j+1/2}}{S_{j+1}^n}
&= \bigg(
- 
\frac{A_{j+1/2}X_{j+1}^n}
{\rho_X-X_{j+1}^n}
\biggl(
\vhs(X_{j+1}^n) - \frac{D(X_{j+1}^n)}{h}
\biggr)
- \Qe^n
\bigg)^-
\\%%%%%%%%
&\geq 
 - A_{j+1/2}\hat{\rho}
\hat{X}\vhs(0)
- \Qe^n.
\end{align*}
Following the same procedure as in the proof of Lemma~\ref{lemmaC} (with $\boldsymbol{\Upsilon}_{\mathrm{e},j}^{\bS,n}$ is defined analogously to~\eqref{eq:Upsilon_e}), we now get
\begin{align*}
\pp{{\Upsilon}_{\mathrm{e},j}^{\bS,n}}
{S_{j+1}^n} 
&= 
 \frac{A_{j+1}}{A_{j}}    + \lambda_{j}\pp{\Psi_{\mathrm{e},j}^{\bS,n}}{S_{j+1}^n}
\geq  \frac{A_{j+1}}{A_{j}}  
- \lambda_{j}\biggl\{
A_{j+1/2}\hat{\rho}_X\hat{X}\vhs(0) + \Qe^n
\\
&\quad 
+ \hat{\rho}\biggl((\rho_X+\hat{X})\Qu^n
+ A_{j+3/2}\biggl(\vhs(0) + \frac{D(\hat{X})}{h}\biggr)\hat{X}
\biggr)
\biggr\}
\\
&\geq \frac{A_{j+1}}{A_{j}}\biggl(
1 - { \tau }
 \frac{{M_1} (\hat{\rho}_X\hat{X}\vhs(0) + \nu_3 )}{h}
\biggr) 
\geq 0. 
\end{align*}
For the update of the effluent concentrations, we get the same result for $\bS_\mathrm{e}^n$ as for~$\bC_\mathrm{e}^n$.

The proof of positivity can be done as in the proof of Lemma~\ref{lemmaC}.
\end{proof}

\begin{lemma}\label{lem:ODEstepC}
If $\boldsymbol{\check{C}}_j^{n+1},\boldsymbol{\check{S}}_j^{n+1}\geq\bzero$ for all $j$ and \eqref{cfl} holds, then the second step of the splitting scheme~\eqref{eq:ODEstep}, \eqref{eq:ODEstep_f} and \eqref{eq:ODEstep_e} (and the analogous formulas for $\bS_j^{n+1}$) satisfy $\bC_j^{n+1},\bS_j^{n+1}\geq\bzero$ for all~$j$.
\end{lemma}

\begin{proof}
Component $k$ of~\eqref{eq:ODEstep} can be written and estimated by means of~\eqref{eq:reactionrates} and \eqref{eq:positivity} as
\begin{align*}
C_j^{(k),n+1}&={\check{C}}_j^{(k),n+1} + 
\tau\sum_{l\in I_{\bC,k}^+}\sigma_{\bC}^{(k,l)}r^{(l)} (\boldsymbol{\check{\bC}}_j^{n+1}, \boldsymbol{\check{\bS}}_j^{n+1})
+\tau\sum_{l\in I_{\bC,k}^-}\sigma_{\bC}^{(k,l)}\bar{r}^{(l)} (\boldsymbol{\check{\bC}}_j^{n+1}, \boldsymbol{\check{\bS}}_j^{n+1}){\check{C}}_j^{(k),n+1}\\
&\geq{\check{C}}_j^{(k),n+1}\Biggl(
1+\tau\sum_{l\in I_{\bC,k}^-}\sigma_{\bC}^{(k,l)}\bar{r}^{(l)} (\boldsymbol{\check{\bC}}_j^{n+1}, \boldsymbol{\check{\bS}}_j^{n+1})
\Biggr)
\geq{\check{C}}_j^{(k),n+1}(1-\tau M_{\bC})\geq 0,
\end{align*}
where the latter inequality follows from~\eqref{cfl}.
Consider now the formula~\eqref{eq:ODEstep_f} for the fill case~(a) with the non-zero coefficient $\omega_{\bar{j}^n+1}^{\mathrm{f},n} :=  A_{\bar{j}^{n}} / A_{\bar{j}^{n}+1}$ (the others are trivial).
Then component $k$ of~\eqref{eq:ODEstep_f} is estimated by using the rewriting above applied to both reaction terms:
\begin{align*}
C_{\bar{j}^n+1}^{(k),n+1}&={\check{C}}_{\bar{j}^n+1}^{(k),n+1}
+ \tau\biggl(\alpha^n \frac{A_{\bar{j}^{n}}}{A_{\bar{j}^{n}+1}} R^{(k)}_{\bC}\biggl(\frac{\check{\bC}_{\bar{j}^n}^{n+1}}{\alpha^n}, \frac{\check{\bS}_{\bar{j}^n}^{n+1}}{\alpha^n}\biggr)
+ R^{(k)}_{\bC}\bigl(\check{\bC}_{\bar{j}^n+1}^{n+1}, \check{\bS}_{\bar{j}^n+1}^{n+1}\bigr)\biggr)\\
&\geq {\check{C}}_{\bar{j}^n+1}^{(k),n+1}\left(
1-\tau\left(\frac{A_{\bar{j}^{n}}}{A_{\bar{j}^{n}+1}} + 1
\right)M_{\bC}\right)
\geq {\check{C}}_{\bar{j}^n+1}^{(k),n+1}\big(
1-\tau (M_3+1)M_{\bC}\big)\geq 0,
\end{align*}
where the latter inequality follows from~\eqref{cfl}.
For fill case~(b), we get in the case $\omega_{\bar{j}^n+1}^{\mathrm{f},n} := \eta^{n+1}$
\begin{align*}
C_{\bar{j}^n+1}^{(k),n+1}
&\geq {\check{C}}_{\bar{j}^n+1}^{(k),n+1}\left(
1-\tau\eta^{n+1}\left(1 + \frac{A_{\bar{j}^{n}+1}}{A_{\bar{j}^{n}}}
\right)M_{\bC}\right)\\
&\geq {\check{C}}_{\bar{j}^n+1}^{(k),n+1}\left(
1-\tau \left(\frac{A_{\bar{j}^{n}}+A_{\bar{j}^{n}+1}}{\alpha^{n+1}A_{\bar{j}^{n}}+A_{\bar{j}^{n}+1}}
\right)M_{\bC}\right)\geq {\check{C}}_{\bar{j}^n+1}^{(k),n+1}\big(
1-\tau (M_3+1)M_{\bC}\big)\geq 0.
\end{align*}
The other cases are similar, for example, in case~(c) with $\omega_{\bar{j}^n+1}^{\mathrm{f},n} := \theta^{n+1}$ one gets
\begin{align*}
C_{\bar{j}^n+1}^{(k),n+1}
&\geq {\check{C}}_{\bar{j}^n+1}^{(k),n+1}\left(
1-\tau\left(\frac{A_{\bar{j}^{n}}+A_{\bar{j}^{n}+1}} {\alpha^{n+1}A_{\bar{j}^{n}-1}+A_{\bar{j}^{n}}+A_{\bar{j}^{n}+1}}
\right)M_{\bC}\right)
\geq {\check{C}}_{\bar{j}^n+1}^{(k),n+1}\left(
1-\tau M_{\bC}\right)\geq 0.
\end{align*}
The cases of extraction are the same and so are the analogous ones for~$\bS_{j+1}^n$.
\end{proof}

\begin{lemma}\label{lem:ODEstepXhat}
Consider the second step of the splitting scheme with formula~\eqref{eq:ODEstep} and the formula obtained by summing all vector components.
If $\smash{\check{X}_j^{n+1}:=\check{C}^{(1),n+1}+\cdots+\check{C}^{(k_{\bC}),n+1} \leq\hat{X}}$ and \eqref{cfl} holds, then the second step of the splitting scheme corresponding to~\eqref{eq:ODEstep} satisfies
\begin{equation}\label{eq:inlemma}
X_j^{n+1}=\check{X}_j^{n+1} + \tau \tilde{R}_{\bC} \bigl(\boldsymbol{\check{\bC}}_j^{n+1}, \boldsymbol{\check{\bS}}_j^{n+1} \bigr)
\leq\check{X}.
\end{equation}
The analogous statement is true for the other update formulas near the surface.
\end{lemma}

\begin{proof}
If $\check{X}_j^{n+1}\geq\hat{X}-\varepsilon$, then $\tilde{R}_{\bC} (\boldsymbol{\check{\bC}}_j^{n+1}, \boldsymbol{\check{\bS}}_j^{n+1})=0$ by assumption~\eqref{eq:techRC}.
Otherwise $\check{X}_j^{n+1}<\hat{X}-\varepsilon$ and \eqref{cfl} gives that $\smash{X_j^{n+1}}$  is estimated by
\begin{equation*}
X_j^{n+1}<\hat{X}-\varepsilon+\tau\tilde{M}\leq\hat{X}.
\end{equation*}
The rest of the update formulas are treated in the same way and the resulting formula contains an additional factor $\omega^{\mathrm{f},n}_j\in[0,1]$ multiplying  $\tau$, so the result also holds for that formula.
\end{proof}

\subsection{Numerics for full mixing}\label{sec:numfullmix}

During the react stage, mixing occurs due to aeration or the movement of an impeller; see Figure~\ref{fig:SBRcycle}.
Then there is no relative velocity between the solids and the liquid, but reactions take place, and the time-dependent concentrations for the mixture below the surface are given by the ODEs~\eqref{eq:ODEmodel}.

Suppose the (PDE or numerical) solution $\bC(z,T_0)$ is known at $t=T_0=t_{n_0}$ when a period of complete mixing starts.
The initial concentrations for the ODEs~\eqref{eq:ODEmodel} are defined as the averages (for $k=1,\dots,k_{\bC}$; analogously for $\bS$)
\begin{equation*}
C^{(k)}(T_0):= \frac{1}{\bar{V}(T_0)}\int_{\bar{z}(T_0)}^{B}A(\xi)C^{(k)}(\xi,T_0)\,\rmd\xi
\approx \frac{h}{\bar{V}(T_0)}\sum_{i=1}^{N}A_{i}C_i^{(k),n_0}
=:C_{\rm aver}^{(k),n_0}.
\end{equation*}
The ODE system~\eqref{eq:ODEmodel} can then be time integrated.
If an ODE mixing period ends at $t=t_{\tilde{n}}$ and the PDE model is to be simulated thereafter, then the total mass below the surface is distributed among the cells by
\begin{equation*}
C_j^{(k),\tilde{n}}:=
\begin{cases}
0 &\text{for $j=1,\ldots,\bar{j}^{\tilde{n}}-1$,} \\
\alpha^{\tilde{n}}C_{\rm aver}^{(k),\tilde{n}}
=(z_{\bar{j}^{\tilde{n}}+1/2}-\bar{z}^{\tilde{n}})C_{\rm aver}^{(k),\tilde{n}} 
\quad  & \text{for $j=\bar{j}^{\tilde{n}}$,}\\
C_{\rm aver}^{(k),\tilde{n}}
\quad  & \text{for $j=\bar{j}^{\tilde{n}}+1,\dots,N$.} 
\end{cases}
\end{equation*}
Recall that the ODEs~\eqref{eq:ODEmodel} are averages of the PDE model~\eqref{finalmod_a} and \eqref{finalmod_b} when the convective and diffusive terms are zero \cite{bcdp_part1}.
Therefore, if time integration is made with the explicit Euler method and~\eqref{cfl} holds, then also the numerical approximations of the ODE solutions to~\eqref{eq:ODEmodel} belong to $\Omega$.

\section{Numerical simulations} \label{sec:numex}

\begin{table}[t]
\caption{List of ASM1 variables of the biokinetic reaction model.}
\label{table:AMS1_vari}
\begin{center}
\begin{tabular}{lcc} \toprule
 \multicolumn{1}{c} { Material } & Notation & Unit \\
\midrule
Particulate inert organic matter 			   & $X_{\rm I}$  & $\rm (g \ COD)\,m^{-3}$\\
 Slowly biodegradable substrate 				   & $X_{\rm S}$  & $\rm (g \ COD)\,m^{-3}$\\
 Active heterotrophic biomass					   & $X_{\rm B, H}$& $\rm (g \ COD)\,m^{-3}$\\
 Active autotrophic biomass 					   & $X_{\rm B, A}$& $\rm (g \ COD)\,m^{-3}$\\
 Particulate products arising from biomass decay & $X_{\rm P}$  & $\rm (g \ COD)\,m^{-3}$\\
 Particulate biodegradable organic nitrogen      & $X_{\rm ND}$ & $\rm (g \ N)\,m^{-3}$\\
 Soluble inert organic matter 				   & $S_{\rm I}$  & $\rm (g \ COD)\,m^{-3}$\\
 Readily biodegradable substrate 				   & $S_{\rm S}$  & $\rm (g \ COD)\,m^{-3}$\\
 Oxygen 										   & $S_{\rm O}$  & $\rm -(g\ COD)\,m^{-3}$\\
 Nitrate and nitrite nitrogen 				   & $S_{\rm NO}$ & $\rm (g \ COD)\,m^{-3}$\\
 $\mathrm{NH}_{4}^{+}+\mathrm{NH}_{3}$ nitrogen  & $S_{\rm NH}$ & $\rm (g \ COD)\,m^{-3}$\\
 Soluble biodegradable organic nitrogen 		   & $S_{\rm ND}$ & $\rm (g \ COD)\,m^{-3}$\\
\bottomrule 
\end{tabular}
\end{center}
\end{table}%

The numerical method in Section~\ref{subsec:scheme} with condition~\eqref{cfl} is first used for the simulation of an SBR process.
In \cite{bcdp_part1}, we demonstrated the process by letting the reaction terms model a denitrification process, which occurs when there is no oxygen present.
Here, the reactions are a modified ASM1 model without alkalinity; see Table~\ref{table:AMS1_vari} for the six particulate and six soluble state variables, and Appendix~A for the reaction terms.
Furthermore, we assume that the mixing during the react stage is achieved  by aeration, which means that the liquid is saturated with dissolved oxygen; a concentration about $S_{\rm O}=10\,\rm g/m^3$.
The constitutive functions used for sedimentation and compression are 
\begin{equation*}
\vhs(X) := \frac{v_0}{ 1 + (X/ \breve{X})^\eta}, \qquad
\sigma_{\mathrm{e}}(X):=\alpha \chi_{\{ X \geq X_{\mathrm{c}} \}} (X-X_{\mathrm{c}}),
\end{equation*}
with $v_0 = 1.76\times 10^{-3}\, \rm m/s$, $\breve{X} = 3.87\, \rm kg/m^3$, $\eta = 3.58$, $X_{\mathrm{c}} = 5\, \rm kg/m^3$ and $\alpha = 0.2\,\rm m^2/s^2$. 
Other parameters are $\rho_X = 1050\, \rm kg/m^3$, $\rho_L = 998\, \rm kg/m^3$, $g = 9.81 \,\rm m/s^2$, and $B = 3\,\rm m$.

\begin{table}[t]
\caption{Time functions for the simulated SBR cycle. `Model' refers to either PDE~\eqref{finalmod} or ODE~\eqref{eq:ODEmodel}.\label{table:SBRcycle}}
\begin{center} 
{\small 
\begin{tabular}{lcccccc} \toprule 
Stage  & Time period [h]& $\Xf(t)\ [\rm kg/m^3]$ & $\Qf(t)\ [\rm m^3/h]$ & $\Qu(t)\ [\rm m^3/h]$ & $\Qe(t)\ [\rm m^3/h]$ & Model\\
\midrule 
Fill & $0\leq t<1$ & 5 &790 & 0& 0 & PDE\\
React & $1\leq t<3$ & 0 & 0 & 0& 0 & ODE\\
Settle & $3\leq t<5$ & 0 &0 & 0& 0 & PDE\\
Draw & $5\leq t<5.5$ & 0 & 0 & 0& 1570 & PDE\\
Idle & $5.5\leq t<6$ & 0 & 0 & 10& 0 & PDE\\
\bottomrule 
\end{tabular}} \end{center} 
\end{table}%

We simulate one sequence of an SBR with the stages specified in Table~\ref{table:SBRcycle}. 
(This is the same scenario as in the first example of~\cite{bcdp_part1} with a denitrification reaction model and without oxygen supply.)
Condition~\eqref{cfl} implies the step length $\tau= 4.3716\cdot 10^{-5} \,$s.
The initial concentrations have been chosen as
\begin{alignat*}{2} 
      & \bC^0(z) =  \boldsymbol{0}  \,{\rm kg/m^3}   \quad && \text{if  $z<2.0 \,{\rm m}$,} \\
      & \bC^0(z) =   ( 0.8889,\    0.0295,\    1.4503,\    0.0904 ,\   0.7371,\    0.0025 )^{\rm T}    \,{\rm kg/m^3}   \quad && \text{if  $z \geq 2.0 \,{\rm m}$,} \\
      & \bS^0(z) =  \boldsymbol{0}  \,{\rm kg/m^3}   \quad && \text{if  $z<2.0 \,{\rm m}$,} \\
      & \bS^0(z) =   ( 0.0400,\    0.0026,\    0.0,\    0.0333,\    0.0004,\    0.0009 )^{\rm T}    \,{\rm kg/m^3}   \quad && \text{if  $z \geq 2.0 \,{\rm m}$,}       
\end{alignat*}
while the feed concentrations are \cite{Henze1987WR}
\begin{alignat*}{1}
& \bCf(t) = \Xf(t)( 0.1273,\    0.5091,\    0.3055,\    3.1819\cdot 10^{-6},\    0.0,\    0.0582 )^{\rm T}    \,{\rm kg/m^3}, \\
 & \bSf(t) = ( 0.04,\    0.064,\    0.0,\    0.001,\    0.0125   ,\ 0.0101 )^{\rm T}    \,{\rm kg/m^3},
\end{alignat*}
where the total solids feed concentration~$\Xf(t)$ varies with time according to Table~\ref{table:SBRcycle}.

\begin{figure}[tbp!] 
\centering 
\subfigure[Particulate inert organic matter ]{ \includegraphics[width=0.47\textwidth]{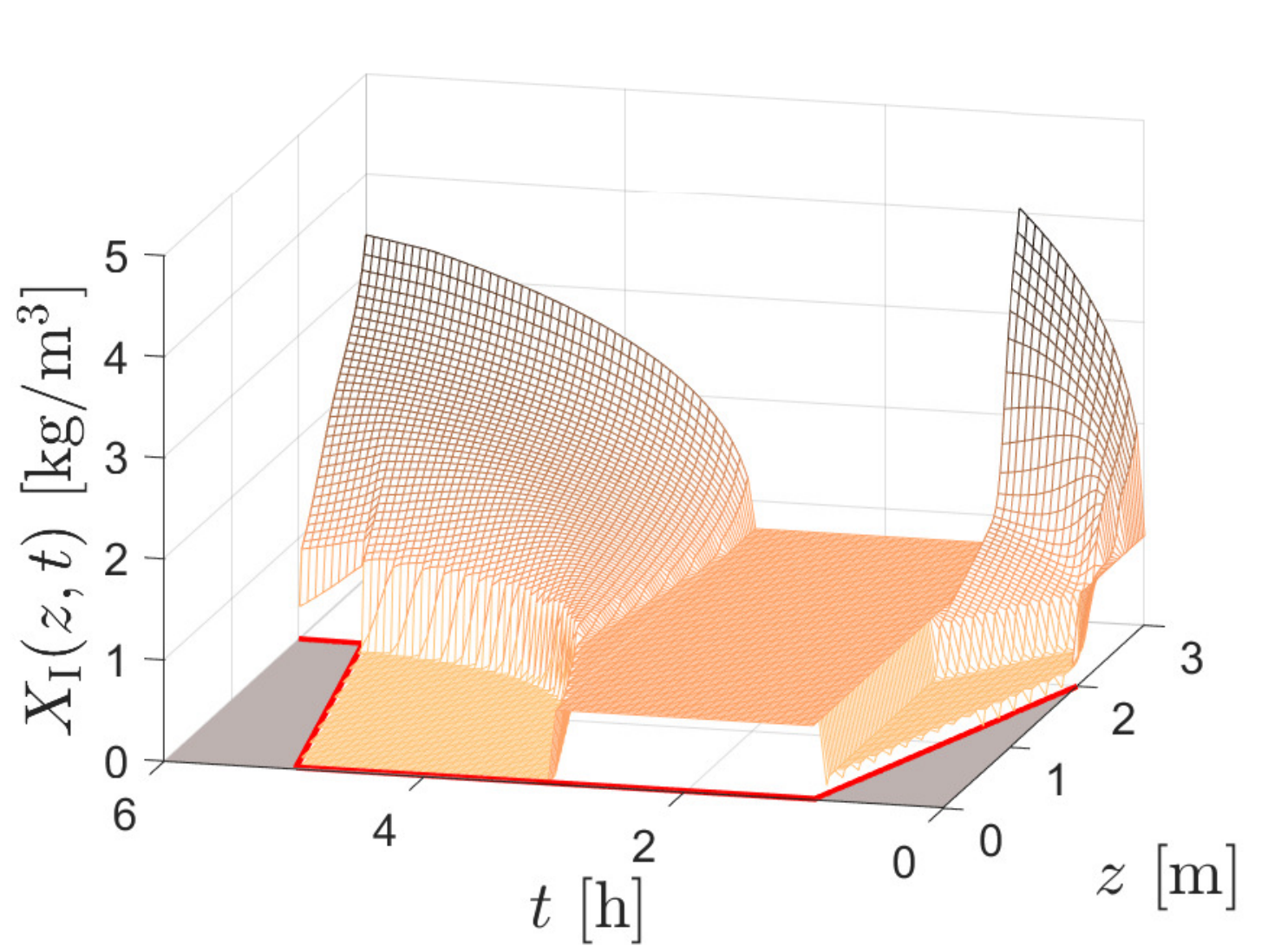}}
\subfigure[Slowly biodegradable substrate]{ \includegraphics[width=0.47\textwidth]{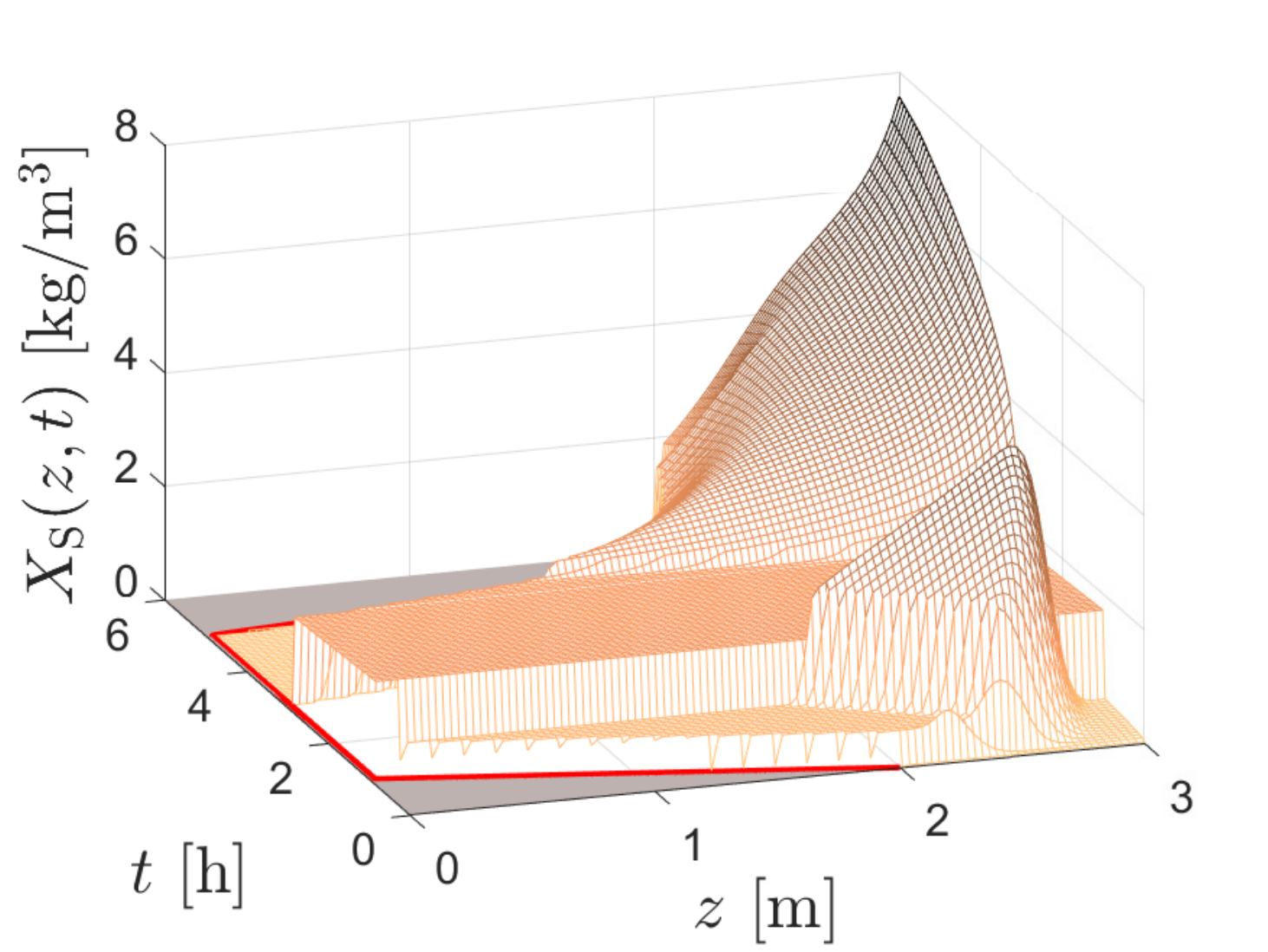}} 
\subfigure[Active heterotrophic biomass ]{ \includegraphics[width=0.47\textwidth]{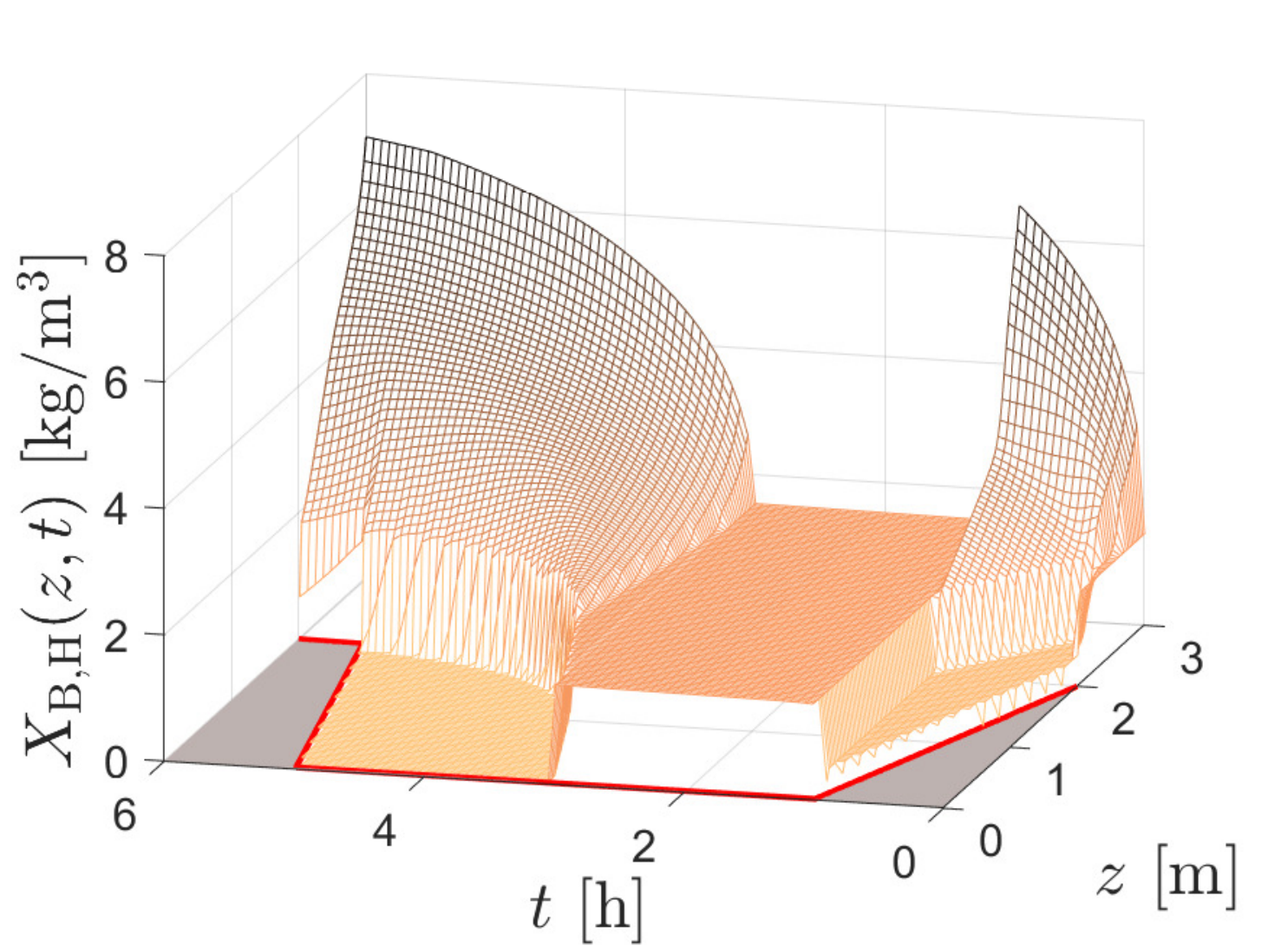}}
\subfigure[Active autotrophic biomass]{ \includegraphics[width=0.47\textwidth]{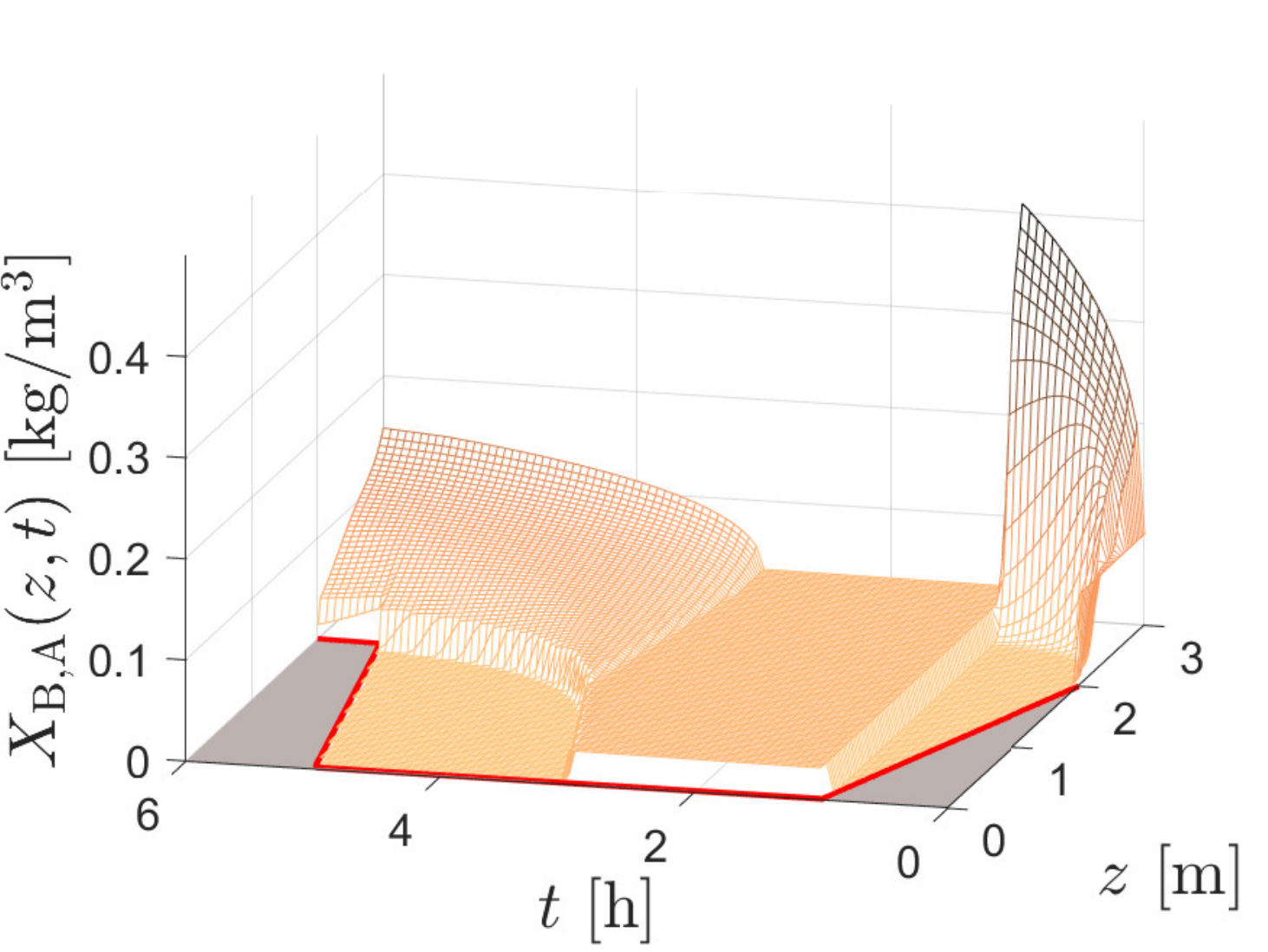}}
\subfigure[Particle products arising from biomass decay]{ \includegraphics[width=0.47\textwidth]{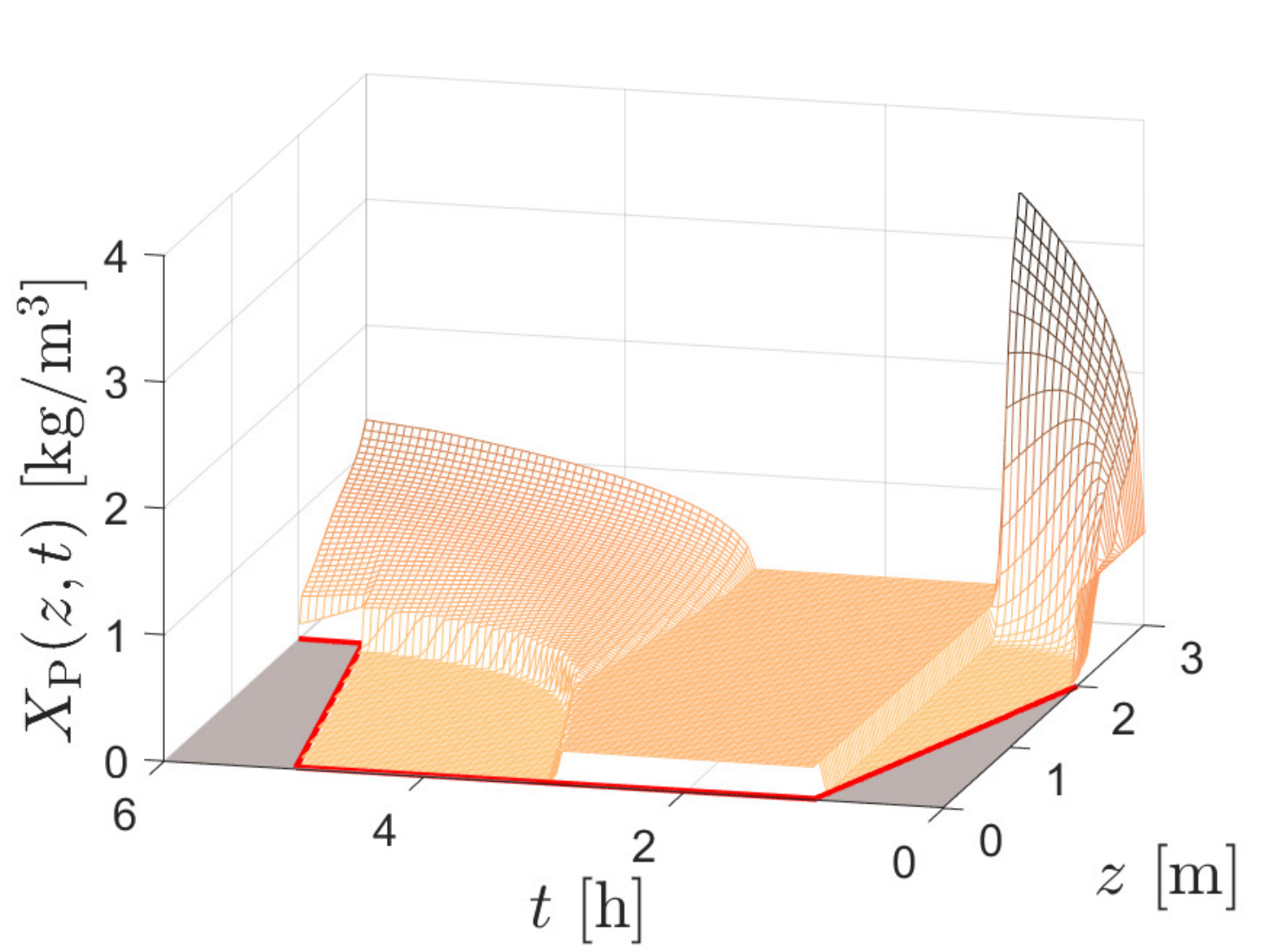}}
\subfigure[Particulate biodegradable organic nitrogen]{ \includegraphics[width=0.47\textwidth]{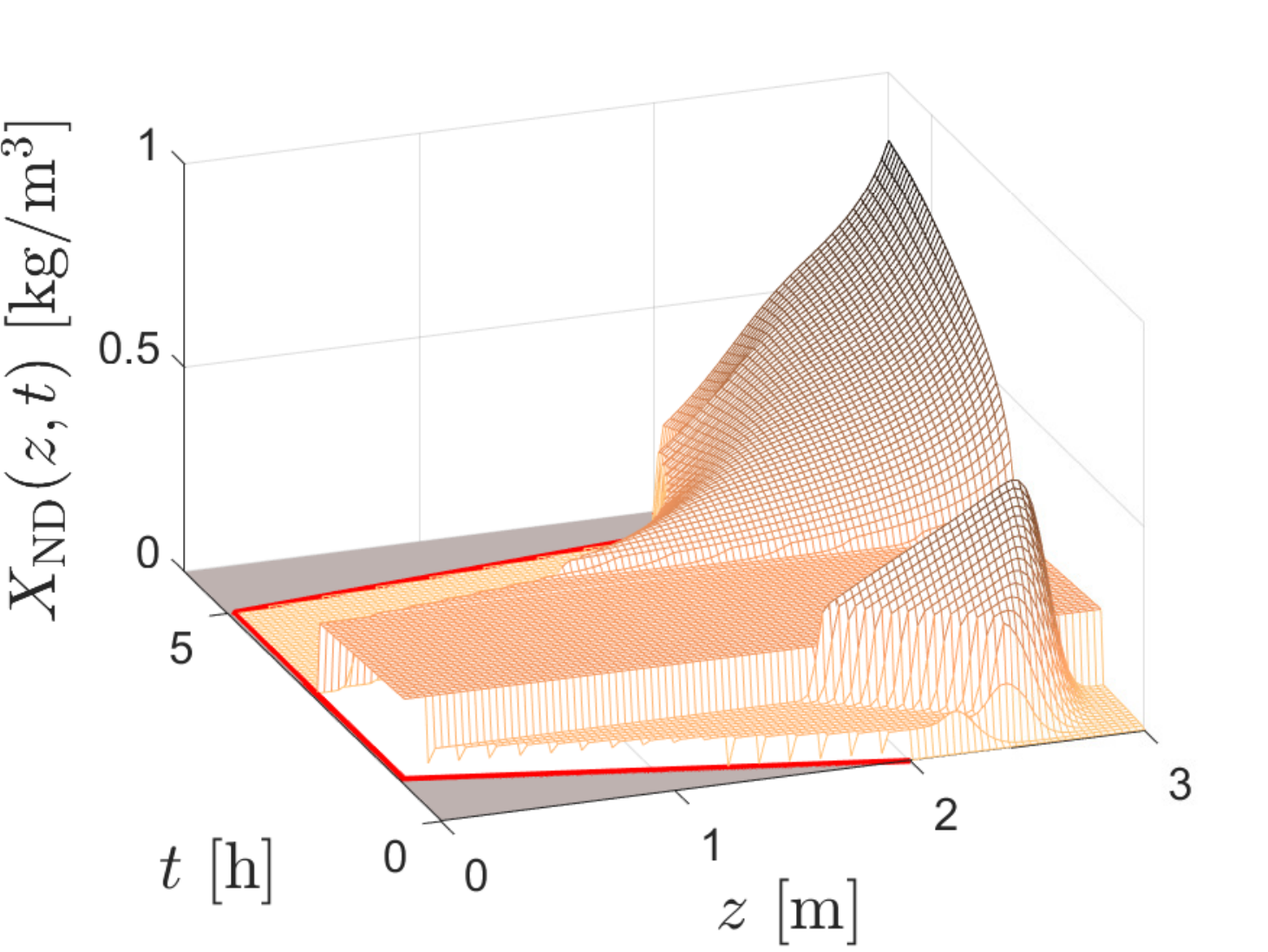}} 
 \caption{Concentrations of the six solid components during a numerical simulation with $N = 100$ until $T = 6\,\rm h$.
For visualization purposes, we do not plot zero  numerical  concentrations above the surface, but fill this region with grey colour.} \label{fig:sim1}
\end{figure}

\begin{figure}[tbp!] 
\centering 
\subfigure[Soluble inert organic matter]{ \includegraphics[width=0.47\textwidth]{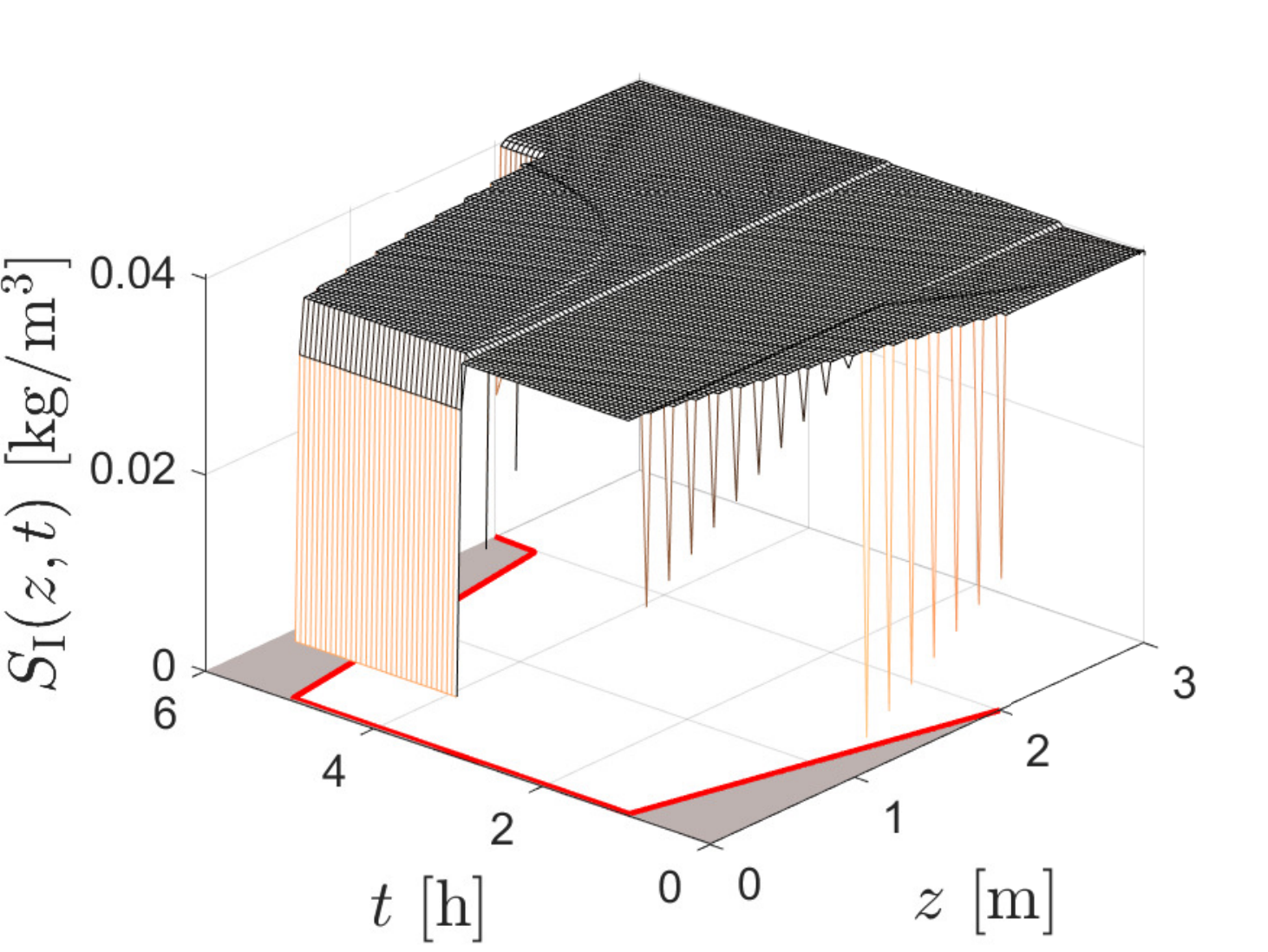}} 
\subfigure[Readily biodegradable substrate]{ \includegraphics[width=0.47\textwidth]{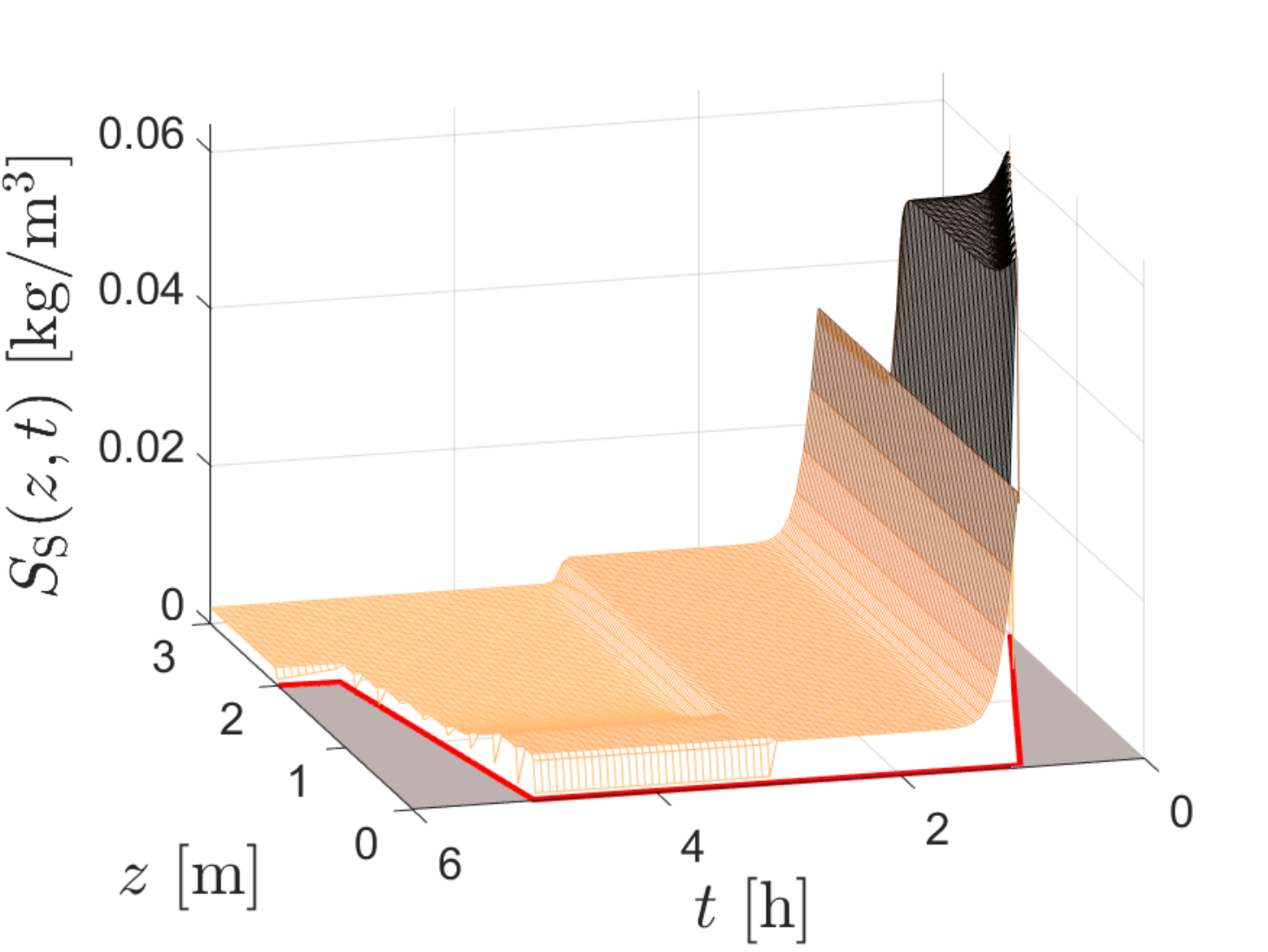}}
\subfigure[Oxygen]{ \includegraphics[width=0.47\textwidth]{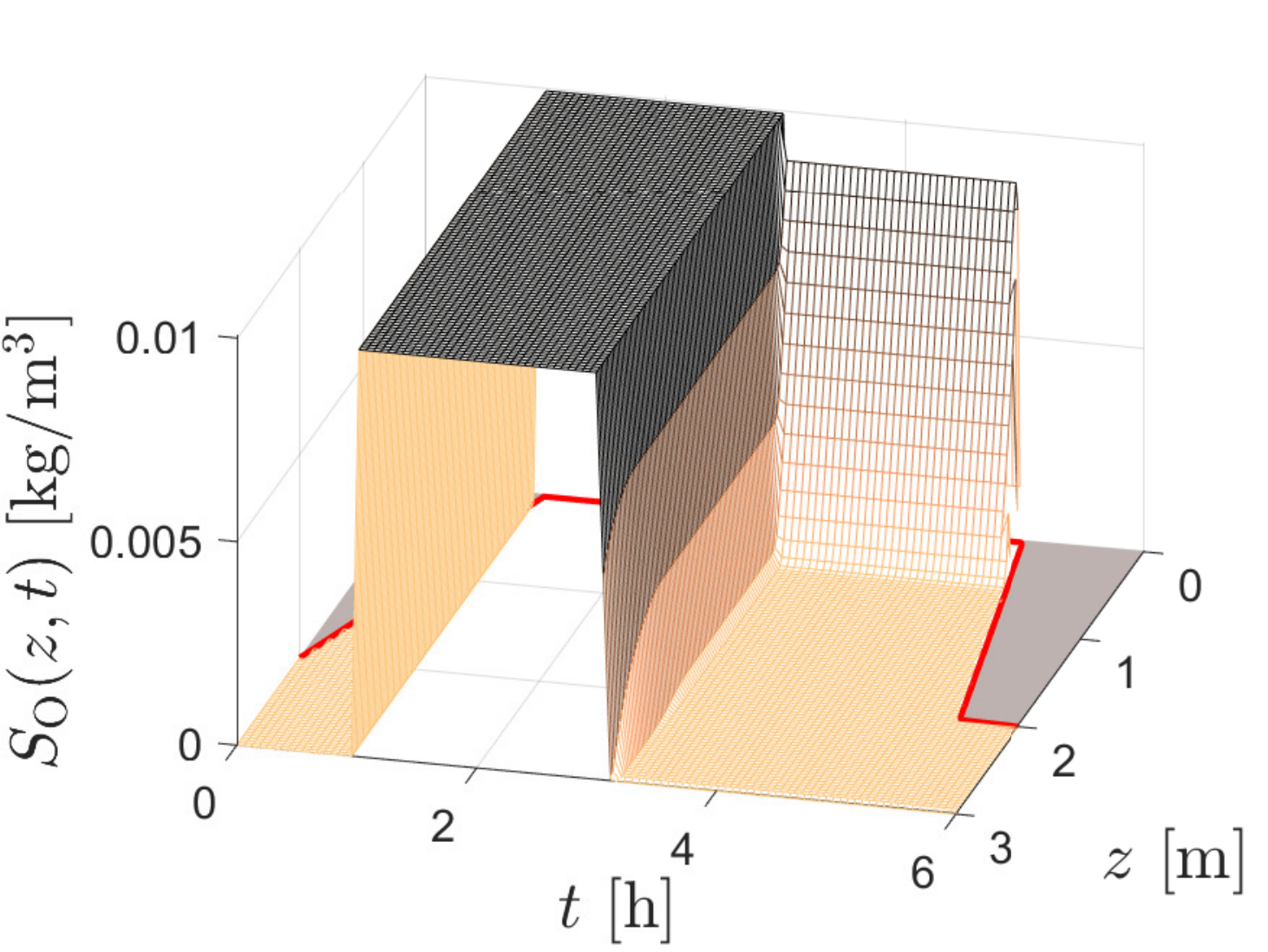}}
\subfigure[Nitrate and nitrite nitrogen]{ \includegraphics[width=0.47\textwidth]{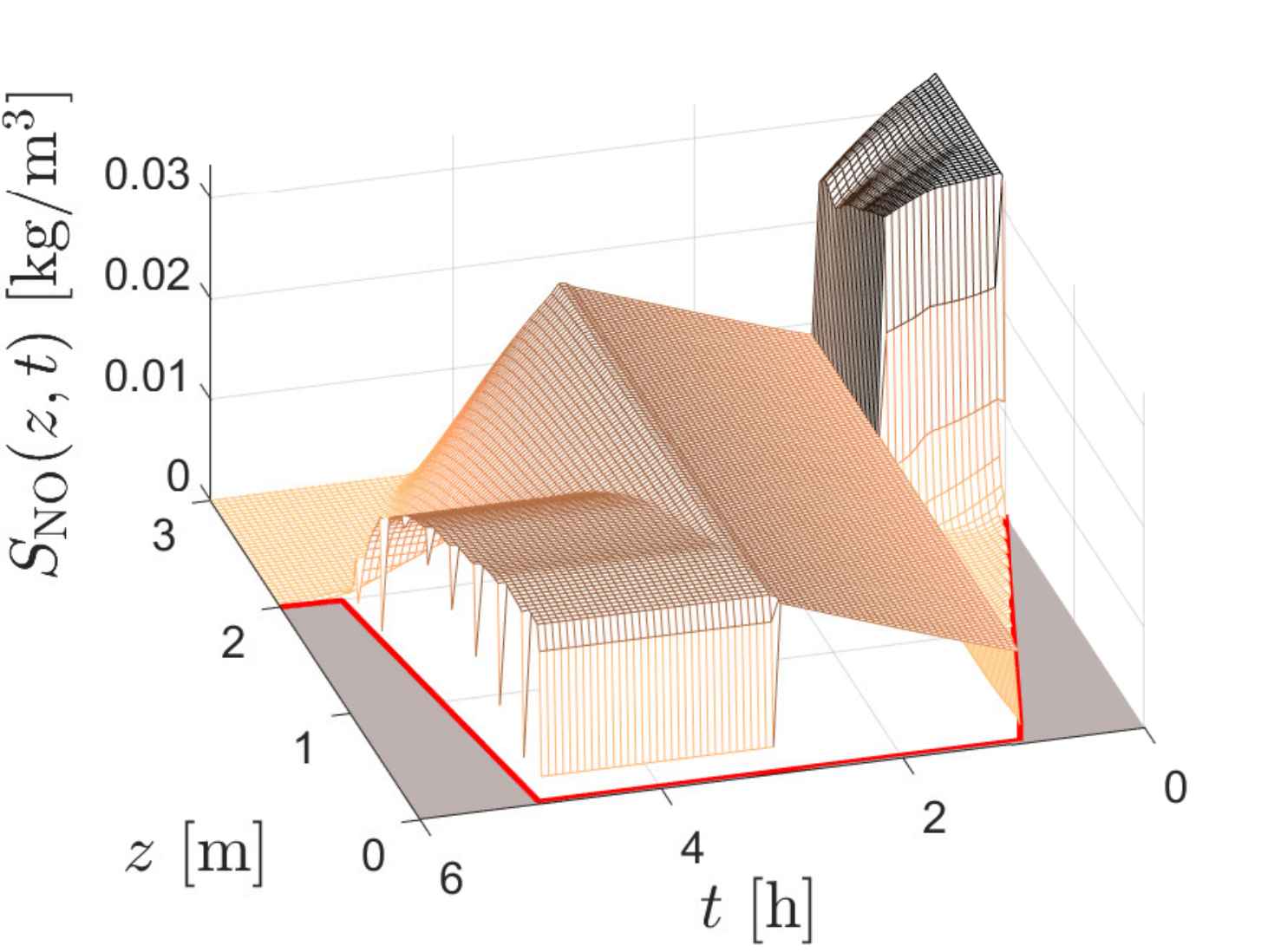}} 
\subfigure[$NH_{4}^+ \ + \ NH_3$ nitrogen]{ \includegraphics[width=0.47\textwidth]{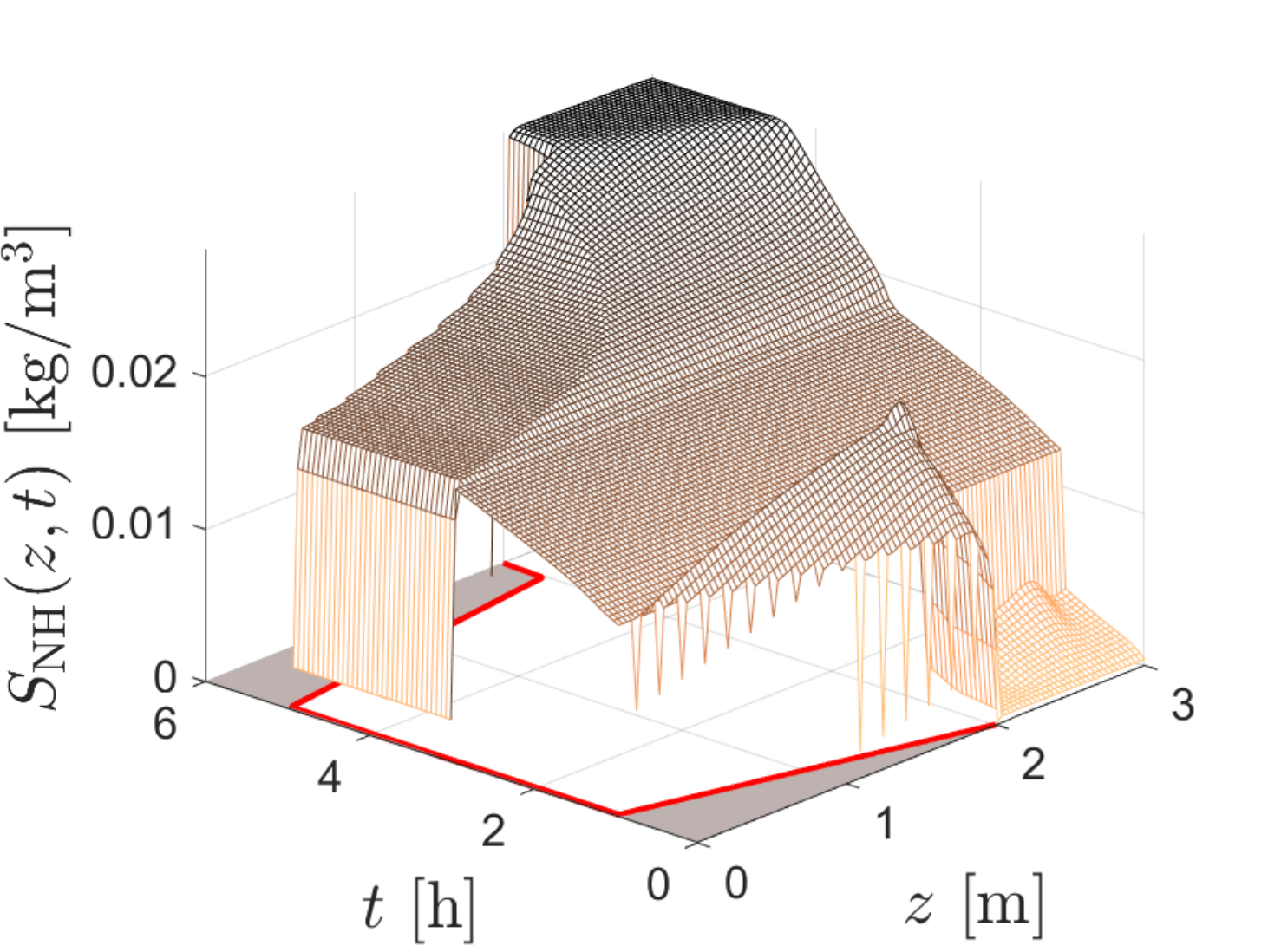}}
\subfigure[Soluble biodegradable organic nitrogen]{ \includegraphics[width=0.47\textwidth]{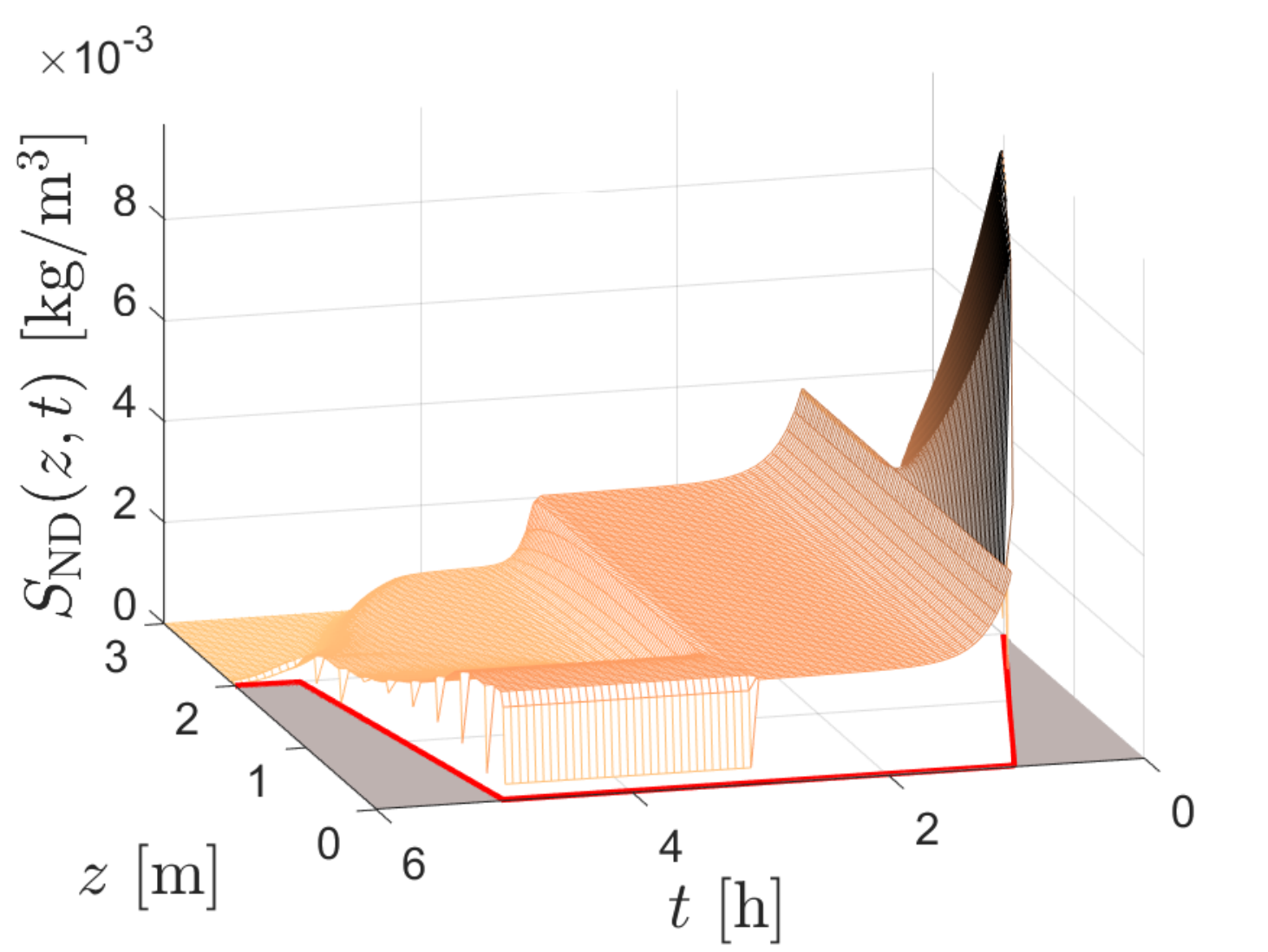}} 
\caption{Concentrations of the six dissolved components during a numerical simulation with $N = 100$ until $T = 6\,\rm h$.
The downwards-pointing peaks at large discontinuities arise because we do not plot zero concentration.} \label{fig:sim2}
\end{figure}

Figures~\ref{fig:sim1} and \ref{fig:sim2} show the simulation results for the concentrations within the vessel of the particulate and soluble components, respectively.
The numerical scheme resolves all the discontinuities accurately.
At time $t=1\,$h, the tank has been filled and the react stage starts and there is  full mixing by aeration.
Consequently, the available components in the tank are distributed homogeneously.
During the react stage, the growth of biomass is slow but the fast consumption of~$S_{\rm S}$ and increase of~$S_{\rm NH}$ are visible  in Figures~\ref{fig:sim2}(b) and (e), respectively.
The nitrification process uses oxygen to produce nitrate and nitrite, which can be seen in plot~(d).
After the react stage, $t\geq 3\,$h, the oxygen is quickly consumed (plot~(c)), but only where there is biomass and $S_{\rm S}$ and~$S_{\rm NH}$ are positive.

\begin{figure}[t] 
\centering 
 \begin{tabular}{ccc}
\includegraphics[scale=0.45]{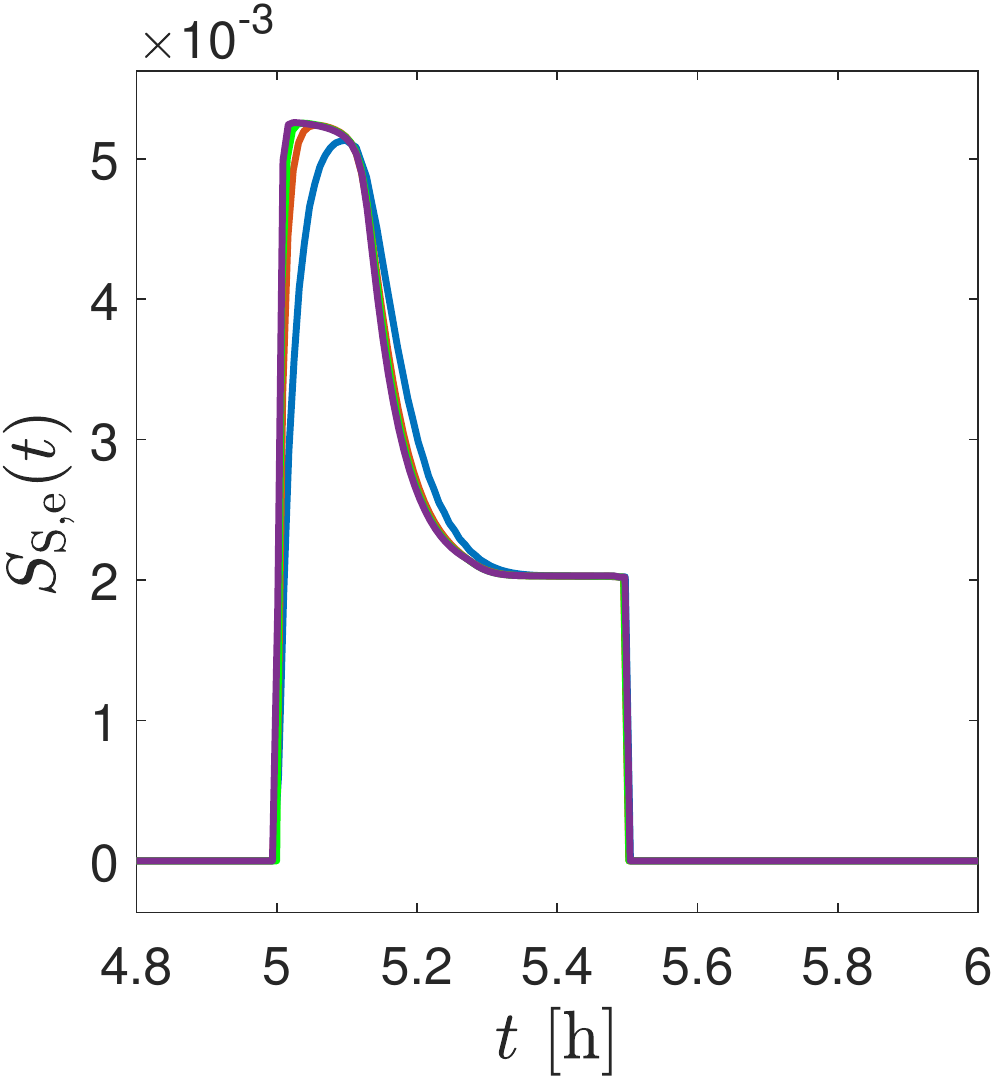} &
\hspace{-0.2cm} \includegraphics[scale=0.45]{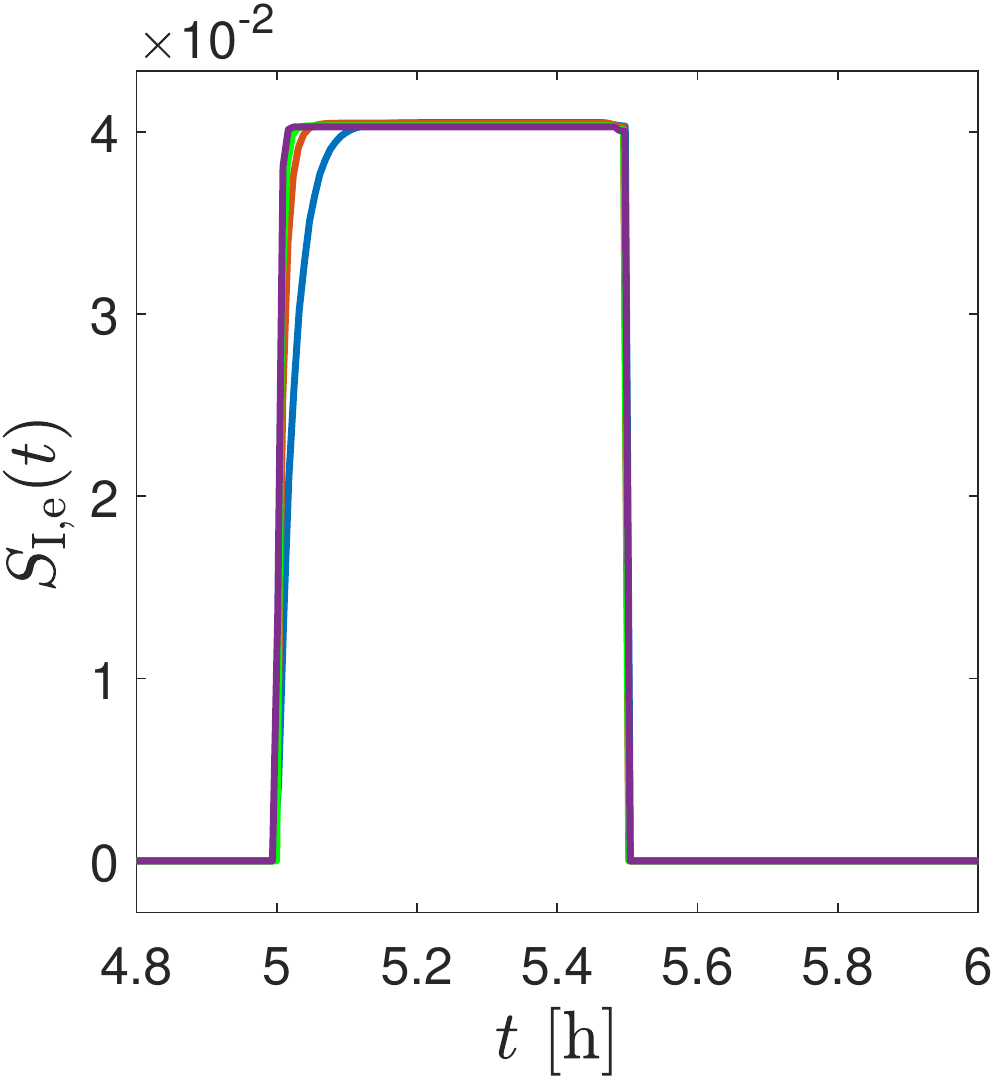} &
\hspace{-0.2cm} \includegraphics[scale=0.45]{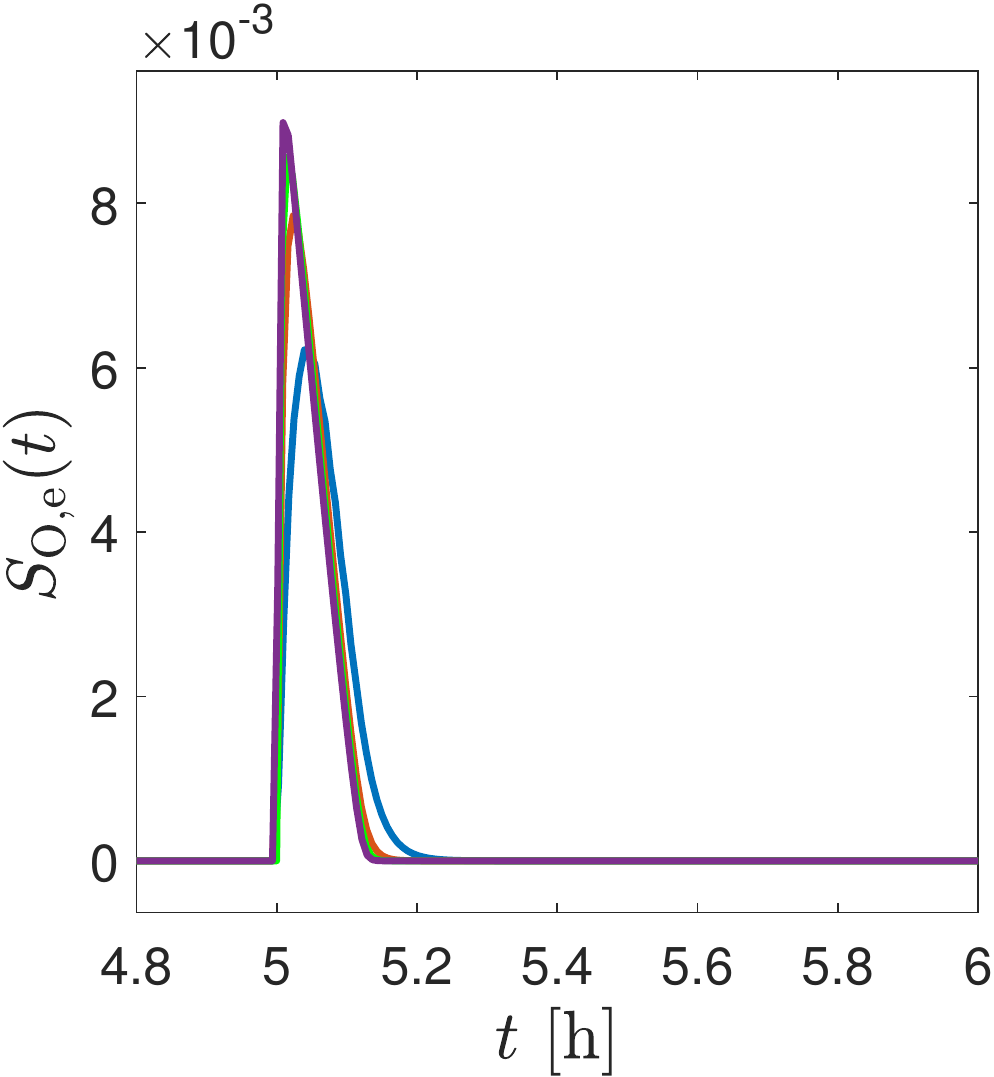} \\
\includegraphics[scale=0.45]{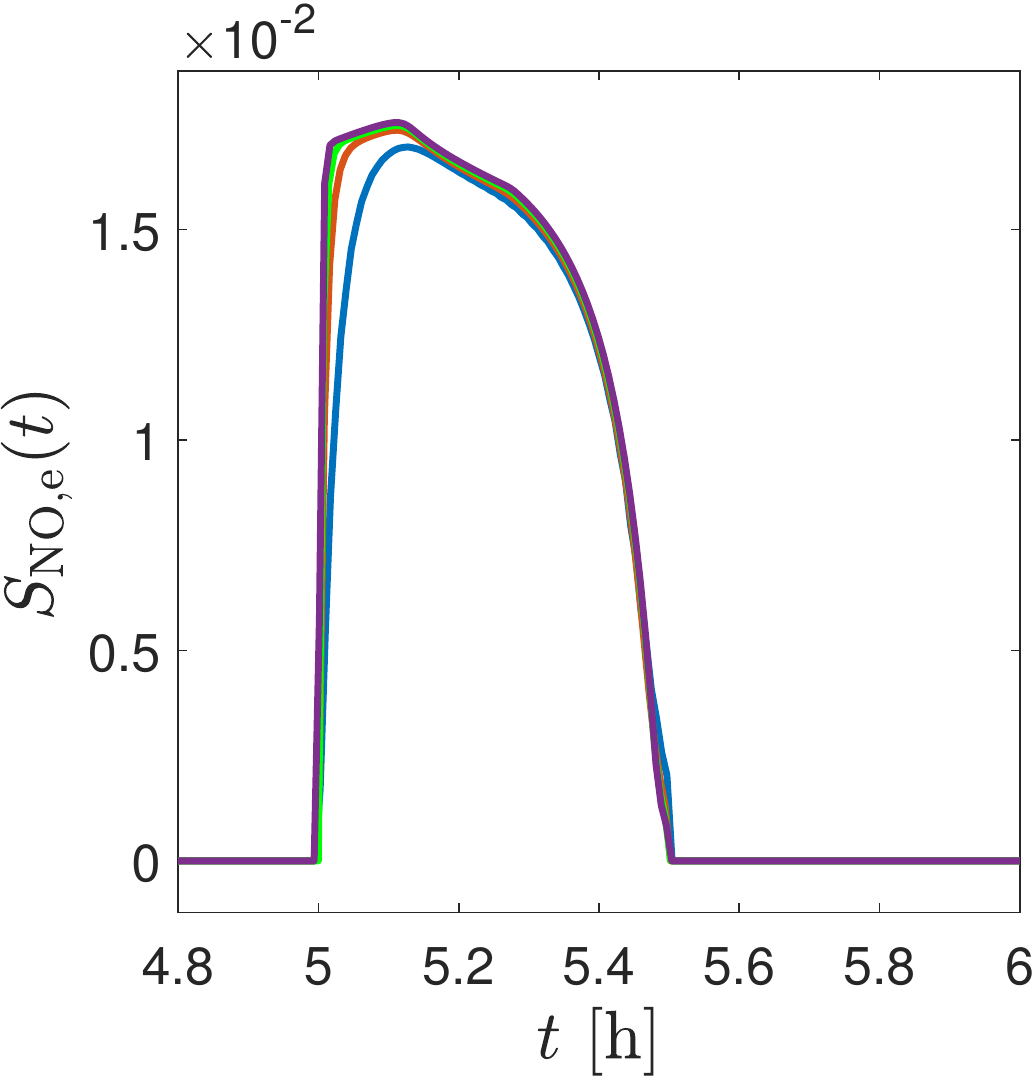} &
\hspace{-0.2cm} \includegraphics[scale=0.45]{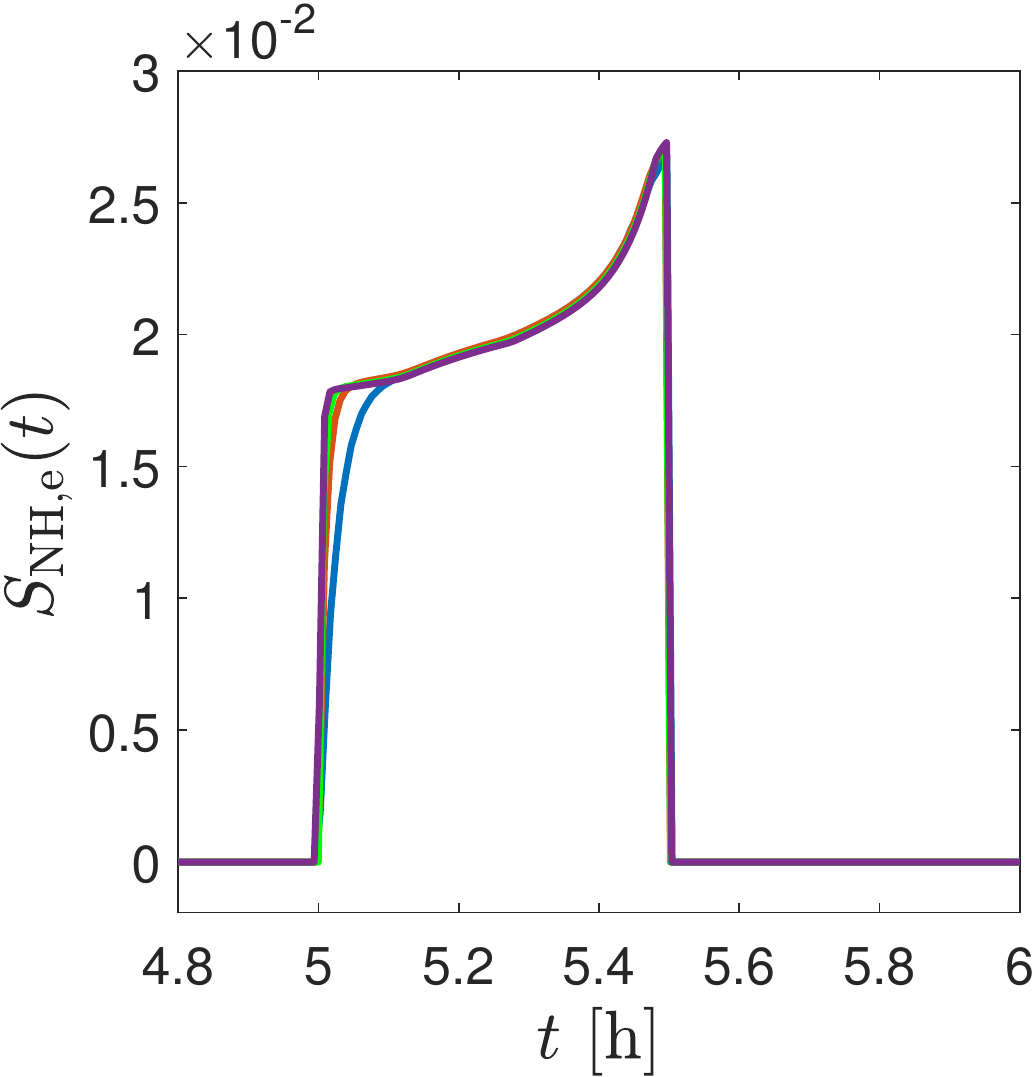} &
\hspace{-0.2cm} \includegraphics[scale=0.45]{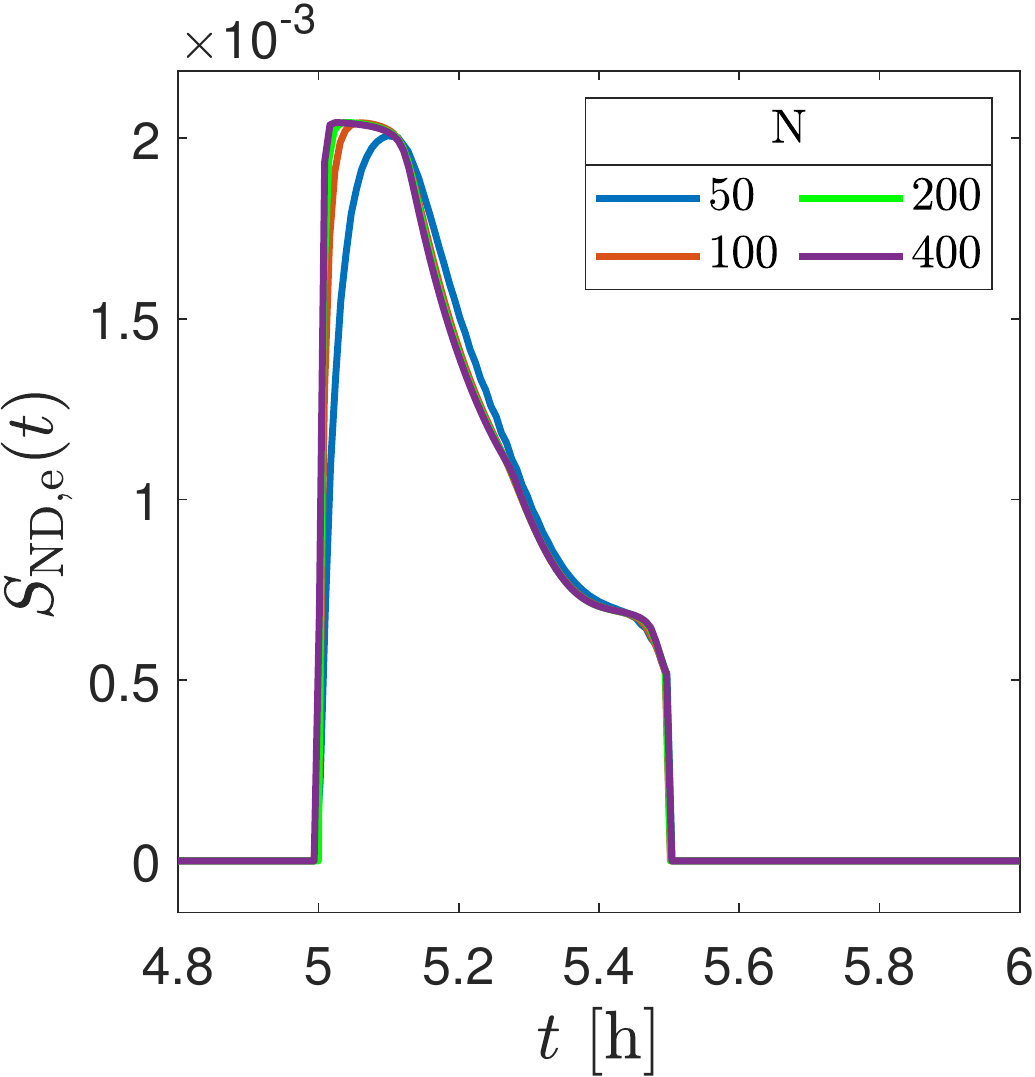} \\
\end{tabular}
 \caption{Simulated effluent concentrations, all in $\mathrm{kg}/\mathrm{m}^3$,  during 
  $t \in [4.8 \, \mathrm{h}, 6 \, \mathrm{h}]$ obtained by the discretizations where the number of computational cells within the tank is $50$, $100$, $200$ and $400$.} \label{fig:effluent}
\end{figure}%

The convergence of the numerical scheme is demonstrated in Figure~\ref{fig:effluent}, where some of the effluent concentrations are plotted.
As was illustrated with the example of  denitrification in~\cite{bcdp_part1} and the present one  
 for the modified  ASM1 model, simulations seem to satisfy the invariant-region property although  we only have a proof of this for the splitting scheme of Section~\ref{subsec:numcfl}.

For the ASM1 example here and a given number of cells~$N$, we calculate the $L^1$ relative difference between the simulation result of the two variants of the numerical scheme according to the following formula, where $\bC_N$ is the result without splitting and~$\bC^{\mathrm{split}}_{N}$ with splitting:
\begin{equation*}
 \mathcal{D}_N(t) := \sum_{k = 1}^{k_{\bC}} \dfrac{\lVert ( C^{(k)}_N-C^{\mathrm{split},(k)}_{N} ) (\cdot,t)\rVert_{L^1(0,B)}}{\lVert C^{\mathrm{split},(k)}_{N} (\cdot,t)\rVert_{L^1(0,B)}}
 +
  \sum_{k = 1}^{k_{\bS}} \dfrac{\lVert ( S^{(k)}_N-S^{\mathrm{split},(k)}_{N} ) (\cdot,t)\rVert_{L^1(0,B)}}{\lVert S^{\mathrm{split},(k)}_{N} (\cdot,t)\rVert_{L^1(0,B)}}.  
\end{equation*}
The result at time $t=T=6\,$h is shown in Table~\ref{table:error}. 
The relative difference is approximately halved as $N$ is doubled.

 \begin{table}[tbp]
\caption{Relative differences~$\mathcal{D}_N(T)$  at the final simulation time $T=6$~h.
\label{table:error}}
\begin{center} 
{\small  
\begin{tabular}{lc} \toprule 
$N$  & $\mathcal{D}_N(T)$\\
\midrule 
50 &  5.8716e-02 \\
100 & 3.5451e-02 \\
200 & 1.6195e-02 \\
400 & 7.5170e-03 \\
800 & 3.3273e-03 \\
\bottomrule 
\end{tabular}} \end{center} 
\end{table}%

\section{Conclusions} \label{sec:concl}

The general model of multi-component reactive settling of flocculated particles given by a quasi-one-dimensional PDE system with moving boundary in Section~\ref{sec:model} was derived in~\cite{bcdp_part1}.
The unknowns are concentrations of biomass particles and soluble substrates, and the reaction terms of the PDE model can be given by any model for  biochemical reactions in wastewater treatment.

The numerical scheme in Section~\ref{subsec:scheme} is designed to ensure conservation of mass across the surface during fill and draw.
Away from the moving boundary, the scheme is the same as in~\cite{SDIMA_MOL}, where it is demonstrated that its order of convergence is not more than one.
The extra treatment near surface will of course not improve that.
An indication of the convergence of the numerical scheme as the mesh size is reduced is demonstrated in Figure~\ref{fig:effluent}. 

The main result of this work is an invariant-region property (Theorem~\ref{theorem}) for the numerical solution if the scheme is computed in a Lie-Trotter-Kato splitting way where the first time step is taken without any reactions and then a step with only the reactions.
Then all the concentrations are nonnegative and the solids concentrations never exceed the maximal packing one.
In particular, the scheme is monotone when the reaction terms are zero.
Simulations with or without splitting have shown to produce very similar outputs and this is demonstrated by the diminishing relative error between such simulations in Table~\ref{table:error}.

A proof of convergence of the method (as $h \to 0$) to a suitably defined weak or entropy weak solution, as well as the corresponding well-posedness (existence and uniqueness) analysis, 
  are still  pending. That said, we point out that available convergence analyses 
   for related strongly degenerate, scalar  PDEs with discontinuous flux (cf., e.g., \cite{Karlsen&R&T2002,Karlsen&R&T2003,Burger&K&T2005a}) rely on the monotonicity 
    of the underlying scheme as well as a uniform bound on the numerical solution, among other properties. 
     Theorem~\ref{theorem} and its proof may be therefore viewed as a partial result to prove  
        convergence of  the numerical scheme presented herein. 
    
Given a moving boundary and a fixed spatial discretization for the numerical scheme, local mass balances have been used to obtain correct update formulas for numerical cells near the surface.
This results in  a scheme with several cases depending on the surface movement.
A certain limitation of the explicit numerical scheme, used without or with splitting, is the restrictive CFL condition (where the time step is esentially proportional to the square of the cell size), implying that very small time steps are needed if accurate approximations on a fine spatial mesh are sought. 
An alternative approach would be to transform the PDE system and have a fixed number of cells below the moving surface.
Such a scheme could possibly also be easier to generalize to a high-order scheme or a more efficient one with semi-implicit time discretization.
The advantage of the present fixed-cell-size numerical scheme is, however, that the model can more easily be generalized to include further sources or sinks at fixed locations, a desirable feature in   applications to wastewater treatment.

\section*{Acknowledgements}
RB~is  supported by ANID (Chile) through 
 projects Centro de Modelamiento Matem\'{a}tico  (BASAL pro\-jects ACE210010 and FB210005); 
  Anillo project ANID/PIA/ACT210030;  CRHIAM, project   ANID/FONDAP/15130015;  and  Fon\-de\-cyt project 1210610. SD~acknowledges support from the Swedish Research Council (Vetenskapsr\aa det, 2019-04601). RP  is supported by ANID scholarship ANID-PCHA/Doctorado Nacional/2020-21200939.

\appendix
\section{The modified ASM1 model}

The alkalinity variable in the original ASM1 model~\cite{Henze2000ASMbook} is removed since it does not influence any other variable.
We introduce an extra Monod factor with a small half-saturation parameter~$\bar{K}_{\rm NH}$ for the concentration $S_{\rm NH}$ in processes nos 1 and 2 (component nos 1 and 2 of the vector $\br$) in order to satisfy condition~\eqref{eq:positivity}- and guarantee non-negative solutions of the ODE system~\eqref{eq:systODE}.
The stoichiometric matrices of the reaction-rate vectors are
\begin{align*}
\bsigmaC &:= \begin{bmatrix}
0&0&0&0&0&0&0&0 \\
0&0&0& 1-f_{\rm P} & 1-f_{\rm P} & 0 & -1 & 0 \\
1 & 1 & 0 & -1 & 0 & 0 & 0 & 0 \\
0 & 0 & 1 & 0 & -1 & 0&0&0 \\
0 & 0 & 0 & f_{\rm P} & f_{\rm P} & 0 & 0 & 0 \\
0 & 0 & 0 & i_{\rm XB}-f_{\rm P}i_{\rm XP} & i_{\rm XB}-f_{\rm P}i_{\rm XP} & 0 & 0 & -1
                  \end{bmatrix},\\
\bsigmaS &:= \begin{bmatrix}
0&0&0&0&0&0&0&0 \\[1mm]
-\dfrac{1}{Y_{\rm H}} & -\dfrac{1}{Y_{\rm H}} &0&0&0&0& 1 & 0\\[4mm]
-\dfrac{1-Y_{\rm H}}{Y_{\rm H}} & 0 & -\dfrac{4.57-Y_{\rm A}}{Y_{\rm A}} &0&0&0&0&0 \\[4mm]
0& -\dfrac{1-Y_{\rm H}}{2.86Y_{\rm H}} & \dfrac{1}{Y_{\rm A}} &0&0&0&0&0 \\[4mm] 
-i_{\rm XB} & -i_{\rm XB} & -i_{\rm XB}-\dfrac{1}{Y_{\rm A}} &0&0& 1 & 0 & 0 \\[1mm]
0&0&0&0&0& -1 & 0 & 1
                  \end{bmatrix},
\end{align*}
and the eight processes are contained in the vector
\begin{equation*}
\br(\bC,\bS):=
\begin{pmatrix}
\mu_{\rm H} \dfrac{S_{\rm N H}}{\bar{K}_{\rm N H}+S_{\rm N H}} \dfrac{S_{\rm S}}{K_{\rm S}+S_{\rm S}} \dfrac{S_{\rm O}}{K_{\rm O, H}+S_{\rm O}} X_{\rm B, H}\\[4mm]
\mu_{\rm H} \dfrac{S_{\rm N H}}{\bar{K}_{\rm N H}+S_{\rm N H}} \dfrac{S_{\rm S}}{K_{\rm S}+S_{\rm S}} \dfrac{K_{\rm O, H}}{K_{\rm O, H}+S_{\rm O}} \dfrac{S_{\rm N O}}{K_{\rm N O}+S_{\rm N O}} \eta_{\rm g} X_{\rm B, H}\\[4mm]
\mu_{\rm A} \dfrac{S_{\rm N H}}{K_{\rm N H}+S_{\rm N H}} \dfrac{S_{\rm O}}{K_{\rm O, A}+S_{\rm O}} X_{\rm B, A}\\[3mm]
b_{\rm H} X_{\rm B, H}\\
b_{\rm A} X_{\rm B, A}\\
k_{\rm a} S_{\rm N D} X_{\rm B, H}\\[3mm]
k_{\rm h} \mu_7(X_{\rm S},X_{\rm B, H})\biggl( \dfrac{S_{\rm O}}{K_{\rm O, H}+S_{\rm O}}+\eta_{\rm h} \dfrac{K_{\rm O, H}}{K_{\rm O, H}+S_{\rm O}} \dfrac{S_{\rm N O}}{K_{\rm N O}+S_{\rm N O}}\biggr)\\[4mm] 
k_{\rm h} \mu_8(X_{\rm B, H},X_{\rm N D})\biggl(  \dfrac{S_{\rm O}}{K_{\rm O, H}+S_{\rm O}}+\eta_{\rm h} \dfrac{K_{\rm O, H}}{K_{\rm O, H}+S_{\rm O}} \dfrac{S_{\rm N O}}{K_{\rm N O}+S_{\rm N O}}\biggr)
\end{pmatrix},
\end{equation*}
where we define
\begin{align*}
\mu_7(X_{\rm S},X_{\rm B, H})
&:=\begin{cases}
0&\text{if $X_{\rm S}=0$ or $X_{\rm B, H}=0$,}\\
 \dfrac{X_{\rm S}X_{\rm B, H}}{K_{\rm X} X_{\rm B,H}+X_{\rm S}}&\text{otherwise,}
\end{cases}
\\
\mu_8(X_{\rm B, H},X_{\rm N D})
&:=\begin{cases}
0 &\text{if $X_{\rm B, H}=0$,}\\
 \dfrac{X_{\rm B, H}X_{\rm N D}}{K_{\rm X}X_{\rm B,H} + X_{\rm S}} &\text{otherwise.}
\end{cases}
\end{align*}
All the constants are given in Table~\ref{table:Name_para}.

\begin{table}[t]
\caption{Stoichiometric and kinetic parameters.}
\label{table:Name_para}
\begin{center} 
\begin{tabular}{cp{8cm}cp{34mm}} \toprule 
 Symbol            & Name & Value & Unit \\ 
\midrule $Y_{\rm A}$   & Yield for autotrophic biomass & 0.24 & $\rm (g\,COD)(g\,N)^{-1}$ \\
 $Y_{\rm H}$   & Yield for heterotrophic biomass & 0.57 & $\rm  (g\,COD)(g\,COD)^{-1}$ \\
 $f_{\rm P}$   & Fraction of biomass leading to particulate products & 0.1 & dimensionless\\
 $i_{\rm XB}$  & Mass of nitrogen per mass of COD in biomass & 0.07 & $\rm (g\,N)(g\,COD)^{-1}$  \\
 $i_{\rm XP}$  & Mass of nitrogen per mass of COD in products from biomass & 0.06 & $\rm (g\,N)(g\,COD)^{-1}$ \\
 $\mu_{\rm H}$ & Maximum specific growth rate for heterotrophic biomass & 4.0 & $\rm d^{-1}$ \\
 $K_{\rm S}$   & Half-saturation coefficient for heterotrophic biomass & 20.0 & $\rm (g\,COD)\,m^{-3}$ \\
 $K_{\rm O, H}$& Oxygen half-saturation coefficient for heterotrophic biomass & 0.25 & $\rm -(g\,COD)\,m^{-3}$ \\
 $K_{\rm NO}$  & Nitrate half-saturation coefficient for denitrifying heterotrophic biomass & 0.5 & $\rm (g\,NO_{3}\text{-}N)\, m^{-3}$ \\
 $b_{\rm H}$   & Decay coefficient for heterotrophic biomass & 0.5 & $\rm d^{-1}$ \\
 $\eta_{\rm g}$& Correction factor for $\mu_{\rm H}$ under anoxic conditions & 0.8 & dimensionless \\
 $\eta_{\rm h}$& Correction factor for hydrolysis under anoxic conditions & 0.35 & dimensionless \\
 $k_{\rm h}$   & Maximum specific hydrolysis rate & 1.5 & $\rm (g\,\text{COD})\, (g\,\text{COD})^{-1}\rm d^{-1}$ \\
 $K_{\rm X}$   & Half-saturation coefficient for hydrolysis of slowly biodegradable substrate & 0.02 & $\rm (g\,\text{COD})(g\,\text{COD})^{-1}$ \\ 
 $\mu_{\rm A}$ & Maximum specific growth rate for autotrophic biomass & 0.879 & $\rm d^{-1}$ \\
$\bar{K}_{\rm NH}$ & Ammonia half-saturation coefficient for aerobic and anaerobic growth of heterotrophs & 0.007 & $\rm (g\,NO_{3}\text{-}N)\, m^{-3}$ \\
 $K_{\rm NH}$ & Ammonia half-saturation coefficient for autotrophic biomass & 1.0 & $\rm (g\,NH_{3}\text{-}N)\, m^{-3}$\\
 $b_{\rm A}$   & Decay coefficient for autotrophic biomass & 0.132 & $\rm d^{-1}$ \\
 $K_{\rm O, A}$& Oxygen half-saturation coefficient for autotrophic biomass & 0.5 & $\rm -(g\,COD)\,m^{-3}$\\
 $k_{\rm a}$   & Ammonification rate & 0.08 & $\rm m^{3}(g COD)^{-1}d^{-1}$  \\
\bottomrule
\end{tabular} \end{center} 
\end{table}

\bibliography{ref_copy}

\end{document}